\documentclass[final,hidelinks,onefignum,onetabnum]{siamart220329}

\usepackage{lipsum,version}
\usepackage{amsfonts}
\usepackage{graphicx}
\usepackage{epstopdf}
\usepackage{algorithmic}
\ifpdf
  \DeclareGraphicsExtensions{.eps,.pdf,.png,.jpg}
\else
  \DeclareGraphicsExtensions{.eps}
\fi

\newsiamremark{remark}{Remark}
\newsiamremark{hypothesis}{Hypothesis}
\crefname{hypothesis}{Hypothesis}{Hypotheses}
\newsiamthm{claim}{Claim}

\headers{Energy stable positivity preserving for Wasserstein gradient flows}{Shiheng Zhang, Jie Shen}

\title{Structure preserving schemes for a class of Wasserstein gradient flows\thanks{Submitted to the editors DATE.
\funding{This work was partially supported by NSFC 12371409.}}}
\author{Shiheng Zhang \thanks{Department of Applied Mathematics, University of Washington, Seattle, WA 98195, USA (\email{shzhang3@uw.edu})}
\and Jie Shen \thanks{School of Mathematical Science,
Eastern Institute of Technology, Ningbo,  China(\email{jshen@eitech.edu.cn})}
}

\usepackage{amsopn}

\usepackage{adjustbox}
\usepackage{graphicx}
\usepackage[numbers]{natbib}
\usepackage{verbatim}
\newcommand{\bx}{\boldsymbol{x}}

\ifpdf
\hypersetup{
  pdftitle=[Structure preserving schemes for Wasserstein gradient flows] {Structure preserving schemes for a class of Wasserstein gradient flows},
  pdfauthor={1}
}
\fi

\begin{document}

\maketitle
\centerline{\em Dedicated to Professor Zhongci Shi}
\medskip

\begin{abstract}
We introduce in this paper  two time discretization schemes tailored for a  range of Wasserstein gradient flows. These schemes are designed to preserve mass, positivity and to be uniquely solvable. In addition, they also ensure energy dissipation in many 
typical scenarios.  Through extensive numerical experiments, we demonstrate the schemes' robustness, accuracy and  efficiency. 
\end{abstract}

\begin{keywords}
Wasserstein gradient flow; positivity preserving;  energy stability; porous media equation
\end{keywords}

\begin{MSCcodes}
65M12; 35K61; 35K55; 65Z05
\end{MSCcodes}

\section{Introduction}
Gradient flows play important roles in mathematical models across various scientific and engineering disciplines, especially in materials science and fluid dynamics \cite{allen1979microscopic, anderson1998diffuse, cahn1958free, doi1988theory, elder2002modeling, gurtin1996two, leslie1979theory, yue2004diffuse}. A general form of gradient flows can be written as
\begin{equation}\label{prb:gradient flow energy}
 \frac{\partial {\rho}}{\partial t}  =  -\mathcal{G} \frac{\delta E}{\delta \rho},
 \end{equation}
 where $E$ is a specific energy functional, and $\mathcal{G}$ is a positive, possibly nonlinear, operator. 
Typical form of $\mathcal{G}$ include: $\mathcal{G} = I$ \cite{allen1979microscopic} or  $\mathcal{G} = -\nabla \cdot M \nabla $ where $M$ is a positive mobility function \cite{cahn1958free}.  
 There are also instances where more intricate metrics are desired. A notable example is the  gradient flows over spaces characterized by the Wasserstein metric, leading to the concept of Wasserstein gradient flows, where $\mathcal{G}$ embodies a nonlinear operator such that 
 \begin{equation}\label{prb:w gradient flow energy}
   \frac{\partial {\rho}}{\partial t}  = \nabla \cdot( \rho  
      \nabla \frac{\delta E}{\delta \rho}).
 \end{equation}
This framework encapsulates a broad class of nonlinear parabolic equations, including the Porous Medium Equation (PME) \cite{vazquez1992introduction,vazquez2007porous}, the Keller-Segel equation \cite{keller1970initiation, keller1971model, patlak1953random}, the Fokker-Planck equation \cite{bolley2011stochastic, sznitman1991topics}, Poisson-Nernst-Planck (PNP) equation \cite{bazant2004diffuse, biler1994debye}, which can all be classified as Wasserstein gradient flows. The surge in interest for Wasserstein gradient flows in recent years is apparent in fields ranging from optimal transport to fluid dynamics and statistical mechanics \cite{villani2021topics,doi2011onsager,peletier2014variational}. Essential properties of the Wasserstein gradient flows include mass conservation, positivity preserving and energy dissipation. 

Designing numerical schemes that can efficiently and accurately capture the essential  properties of the Wasserstein gradient flows is still a challenge. Since the Wasserstein gradient flows are a special class of gradient flows, we can consider established techniques for gradient flows. However, energy stable schemes for gradient flows, such as convex splitting \cite{elliott1993global,eyre1998unconditionally,shen2012second}, IEQ \cite{yang2016linear,zhao2017numerical} and SAV \cite{shen2018convergence,shen2018scalar,shen2019new}), do not inherently guarantee positivity preservation. The celebrated JKO scheme \cite{jordan1998variational} is a widely used approach to solve variational problems associated with  Wasserstein gradient flows. This method's variational structure ensures unconditional energy stability and positivity preserving. However, it requires solving a nonlinear minimization problem in Wasserstein metric. Significant efforts have been made towards  efficient computation of the JKO scheme, such as \cite{carrillo2022primal, carrillo2021variational,  peyre2015entropic,  leclerc2020lagrangian,  jacobs2020fast, jacobs2021back,fu2023high, fu2023highv}. Despite extensive research on the JKO scheme, efficiently solving it remains challenging due to the inherent complexity of the Wasserstein metric. Some other approaches have been developed recently that can preserve positivity and energy dissipation for Wasserstein gradient flows, but they either are not uniquely solvable \cite{bailo2018fully,hu2023positivity}, or are limited to certain   Wasserstein gradient flows  \cite{shen2020unconditionally,ShenXu2021}. We also refer to more recent work \cite{duan2021structure,duan2022convergence} for other efforts in developing structure preserving schemes for Wasserstein gradient flows.

We propose in this paper two approaches to construct time discretization schemes for Wasserstein gradient flows.
\begin{enumerate}
    \item In the first approach, we set \(E(\rho) =  \int_{\Omega} H(\rho) + \rho v \,d\mathbf{x}\), where \(H''(\rho)>0\) ensures convexity, and $v$ is a given potential function. We construct  first-order and second-order schemes that can preserve mass, positivity and are uniquely solvable. Furthermore, the first-order scheme is also energy dissipative in some typical scenarios.
    \item In the second approach, \(E(\rho)\) can take   a more general form, but can be split into two parts: one part is convex, and the other is bounded from below. This assumption is significantly less restrictive than the assumption required by the convex splitting method since only one part is required to be convex. We construct  first-order and second-order schemes by introducing a SAV \cite{shen2019new}, and show that  they preserve mass, positivity, are uniquely solvable and energy dissipative with respect to a modified energy.
\end{enumerate}
These schemes are nonlinear in nature, but their solutions can all be interpreted as  minimizers of strictly convex functionals.

 The paper is structured as follows: We introduce in Section 2  the two  approaches in the semi-discrete form, and establish their essential properties.  We present in  Section 3  a spatial discretization which can inherent essential properties in the space continuous case using the finite difference method as an example. We carry out a series of numerical experiments in Section 4 to validate our schemes's  accuracy and efficiency. Finally, we provide some concluding remarks in  Section 5.

\section{Time Discretization}
In this section, we present two distinct time discretization approaches, each tailored to specific scenarios of the Wasserstein gradient flows \eqref{prb:w gradient flow energy} with a smooth functional $E: \Omega \times \mathbb{R} \rightarrow \mathbb{R}$. We shall assume throughout this paper that the boundary condition for $\rho$  is either periodic or homogeneous Neumann, i.e., $\frac{\partial \rho}{\partial \boldsymbol{n}}=0$, 
and use  $(a, b)$ to denote the integral $\int_{\Omega} ab \mathrm{~d}\bx$.

\subsection{The first approach (S1)}
The Wasserstein gradient flow is as follows:
\begin{equation} \label{wgde1}
    \frac{\partial \rho}{\partial t}= \nabla \cdot\left(\rho \nabla  \frac{\delta E}{\delta \rho} \right), \quad \bx \in \Omega \subset \mathbb{R}^d, t>0.
\end{equation}
Under the assumption that $E(\rho) =  \int_{\Omega} H(\rho) + \rho v \mathrm{~d} \bx$ with $H''(\rho)>0$, and using the identity $\nabla \rho=\rho \nabla \log \rho$, equation \eqref{wgde1} can be reformulated as:
\begin{equation}\label{wgde2}
\frac{\partial \rho}{\partial t} = \nabla \cdot\left(\rho \left( H''(\rho)\rho \nabla \log \rho +\nabla v\right)\right).
\end{equation}

A first-order time-discretized scheme can be written as follows:
 \begin{equation*}\label{wgds1}
        \frac{\rho^{n+1}- \rho^n}{\delta t} = \nabla \cdot\left(\rho^n \left( H''\left(\rho^n\right)\rho^n \nabla \log \rho^{n+1} +\nabla v\right)\right). \tag{S1}
    \end{equation*}
This scheme is a generalization of the schemes introduced in \cite{ShenXu2021} and \cite{gu2020bound}. It is a nonlinear scheme, but can guarantee mass conservation law, positivity, also the energy dissipation law if $H''(\rho)\rho = 1$ or $\nabla v = 0$, as stated below.
\begin{theorem}\label{Thm:S1}
    Assuming  $H''(\rho)>0$ and  $\rho^n > 0$, the  scheme (S1) exhibits the following attributes:
    \begin{itemize}
        \item[1.] Mass conservation:
        $ \int_{\Omega} \rho^n \mathrm{~d} \bx= \int_{\Omega} \rho^{n+1} \mathrm{~d} \bx$;
        \item[2.] Positivity preserving: $\rho^{n+1}>0$;
        \item[3.] Unique solvability;
        \item[4.] Energy dissipation law: if $H(\rho)=\rho(\log \rho-1)+c \rho$, we have
\begin{equation}\label{disc energy diss}
    \begin{aligned}
         &\int_{\Omega} \rho^{n+1} (\log \rho^{n+1}  - 1+ v) \mathrm{~d} \bx-\int_{\Omega} \rho^n (\log \rho^n  - 1+ v) \mathrm{~d} \bx \\&\leq -\delta t \int_{\Omega} \rho^n \left|\nabla (\log \rho^{n+1} + v )\right|^2 \mathrm{~d} \bx;
    \end{aligned}
\end{equation}
 or if $\nabla v = 0$, we have
\begin{equation}\label{disc energy diss2}
    \begin{aligned}
         &\int_{\Omega} \rho^{n+1} (\log \rho^{n+1} - 1) \mathrm{~d} \bx-\int_{\Omega} \rho^n (\log \rho^n - 1) \mathrm{~d} \bx \\&\leq -\delta t \int_{\Omega} \left(\rho^n\right)^2 H''(\rho^n)\left|\nabla \log \rho^{n+1} \right|^2 \mathrm{~d} \bx.
    \end{aligned}
\end{equation}
    \end{itemize}
\end{theorem}
\begin{proof}
    The mass conservation is obtained by integrating equation \eqref{wgds1} over the domain $\Omega$ and utilizing the Neumann boundary conditions or periodic boundary conditions imposed to $\rho^{n+1}$.
    
    The preservation of positivity for $\rho^{n+1}$ is ensured through the inclusion of $\log \rho^{n+1}$ in our formulation. 

    For the unique solvability, we define $\mathcal{L}_n$ as a linear operator by $\mathcal{L}_n g = u$ where $u$ is the solution of the following elliptic equation
    \begin{equation*}
         -\nabla \cdot \left( \left(\rho^n\right)^2H''(\rho^n) \nabla  u   \right) = g, \quad \int_{\Omega} u \mathrm{~d}\bx = 0
    \end{equation*}
    with Neumann boundary conditions or periodic boundary conditions.
    Hence, $\mathcal{L}_n$ is self-adjoint and semi-positive definite. We then consider a functional
    \begin{equation*}
    \begin{aligned}
    F[\rho^{n+1}]&:=\int_{\Omega} \rho^{n+1} (\log \rho^{n+1} -1 ) \mathrm{~d} \bx+\frac{1}{2\delta t} \int_{\Omega}\left(\rho^{n+1}-\rho^{n}\right) \mathcal{L}_n\left(\rho^{n+1}-\rho^{n}\right) \mathrm{~d} \bx\\
          &+\int_{\Omega}\rho^{n+1} \mathcal{L}_n  \nabla \cdot \left(\rho^n\nabla v\right)\mathrm{~d} \bx.
    \end{aligned}
    \end{equation*}
    Notably, $F$ stands as a strictly convex functional within the set $$\mathcal{A}:=\left\{\rho^{n+1} \in H^2(\Omega): \rho^{n+1}>0 \text { in } \Omega \right\}.$$ Remarkably, the Euler-Lagrange equation of $F$ under the constraint of mass conservation is
    \begin{equation}
    \begin{aligned}
        &\log\rho^{n+1} + \frac{1}{\delta t}\mathcal{L}_n\left( \rho^{n+1} - \rho^n\right) + \mathcal{L}_n  \nabla \cdot \left(\rho^n\nabla v\right)= \lambda,\\
        & \int_{\Omega} \rho^{n+1} - \rho^{n} \mathrm{~d}\bx = 0,
    \end{aligned}
    \end{equation} 
    where $\lambda$ is the Lagrange multiplier of the mass conservation. It is equivalent to the equation given in scheme \eqref{wgds1} with the definition of $\mathcal{L}_n$. Consequently, the unique minimizer of $F$ serves as the unique solution to scheme \eqref{wgds1}.

  If $H(\rho)=\rho(\log \rho-1)+c \rho$, we have $H''(\rho)\rho = 1$. Taking the inner product of equation \eqref{wgds1} with $\log \rho^{n+1} + v$ and employing integration by parts, we obtain
    \begin{equation*}
        \int_{\Omega} \left( \rho^{n+1} - \rho^n \right)\left(\log \rho^{n+1} + v\right)\mathrm{~d}\bx =  -\delta t \int_{\Omega} \rho^n \left|\nabla (\log \rho^{n+1} + v )\right|^2 \mathrm{~d} \bx.
    \end{equation*}
To facilitate further simplification, we utilize the following inequality which can be derived by Taylor expansion:
    \begin{equation}\label{eq:log}
        (a-b) \log a=(a \log a-a)-(b \log b-b)+\frac{(a-b)^2}{2 \xi}, \quad \xi \in[\min \{a, b\}, \max \{a, b\}].
    \end{equation}
Subsequently, we can infer:
    \begin{equation*}
        (\rho^{n+1} \log\rho^{n+1}-\rho^{n+1})-(\rho^{n} \log\rho^{n}-\rho^{n}) \leq \left(\rho^{n+1}- \rho^n\right)\log\rho^{n+1}.
    \end{equation*}
Employing this inequality directly leads us to
\begin{equation}
\begin{aligned}
&\int_{\Omega} \rho^{n+1} \left(\log\rho^{n+1}+v - 1\right)-\rho^{n} \left(\log\rho^{n}+v-1\right) \mathrm{~d}\bx\\
&\leq\int_{\Omega} \left( \rho^{n+1} - \rho^n \right)\left(\log \rho^{n+1} + v\right)\mathrm{~d}\bx\\
&=-\delta t \int_{\Omega} \rho^n \left|\nabla (\log \rho^{n+1} + v )\right|^2 \mathrm{~d} \bx\\ &\leq 0,
\end{aligned}
\end{equation}
which is the desired energy dissipation law \eqref{disc energy diss}. 

Similarly,  if $\nabla v= 0$ the energy dissipation law \eqref{disc energy diss2} can be obtained by taking the inner product of equation \eqref{wgds1} with $\log \rho^{n+1}$ and integrating by parts.
\end{proof}
We can construct a second-order scheme by combining the second-order BDF with Adams-Bashforth extrapolation as follows:
\begin{equation}\label{wgds1o2}
        \frac{3\rho^{n+1}-4\rho^n+\rho^{n-1}}{2\delta t} = \nabla \cdot \phi^{*, n+\frac{1}{2}} \nabla \log \rho^{n+1}+\nabla \cdot \rho^{*,n+\frac 12} \nabla v,
    \end{equation}
where $\phi = \rho^2H''(\rho)$, and for any function $\psi$,
\begin{equation}\label{abe}
        \psi^{*, n+\frac{1}{2}}=\begin{cases} 2 \psi^n-\psi^{n-1} & \text{ if } \psi^n\ge \psi^{n-1}\\
\frac {1}{2/\psi^n-1/\psi^{n-1}}& \text{ if } \psi^n< \psi^{n-1}\end{cases},
    \end{equation}
which is a modified Adams-Bashforth extrapolation to preserve positivity, i.e., $\psi^{*, n+\frac{1}{2}} > 0$ if $\psi^{n}, \psi^{n-1} > 0$. To find \(\rho^1\), we can use the the first-order method (S1).
\begin{theorem}\label{thm:s1-2nd}
    Assume that $\rho^1$ is obtained from the first-order scheme. The  second-order scheme \eqref{wgds1o2} exhibits the following attributes:
      \begin{itemize}
        \item[1.] Mass conservation:
        $ \int_{\Omega} \rho^n \mathrm{~d} \bx= \int_{\Omega} \rho^{n+1} \mathrm{~d} \bx$;
        \item[2.] Positivity preserving: $\rho^{n+1}>0$;
        \item[3.] Unique solvability.
    \end{itemize}
\end{theorem}
\begin{proof}
    The proof is similar to that of Theorem \ref{Thm:S1}. We only need to modify slightly the definition of the linear operator $\tilde{\mathcal{L}}_n$ as follows:  $\tilde{\mathcal{L}}_n g = u$ is defined by the following elliptic equation:
    \begin{equation*}
         -\nabla \cdot \left( \phi^{*, n+\frac{1}{2}} \nabla  u   \right) = g, \quad \int_{\Omega} u \mathrm{~d}\bx = 0
    \end{equation*}
     with the homogeneous Neumann boundary conditions or periodic boundary conditions.
     We can also define a slightly different convex functional $\tilde{F}$ as follows:
    \begin{equation*}
    \begin{aligned}
          \tilde{F}[\rho^{n+1}]&:=\int_{\Omega} \rho^{n+1} (\log \rho^{n+1} -1 ) \mathrm{~d} \bx+\int_{\Omega}\rho^{n+1} \tilde{\mathcal{L}}_n  \nabla \cdot \left(\rho^n\nabla v\right)\mathrm{~d} \bx\\&\frac{1}{12\delta t} \int_{\Omega}\left(3\rho^{n+1}-4\rho^{n}+\rho^{n-1}\right) \tilde{\mathcal{L}}_n\left(3\rho^{n+1}-4\rho^{n}+\rho^{n-1}\right) \mathrm{~d} \bx.
    \end{aligned}
    \end{equation*}
    Then it can be shown that the solution of \eqref{wgds1o2} is the unique minimizer of the above convex functional.
\end{proof}

\subsection{The second approach (S2)}
While the scheme S1 retains many desired properties, its applicability is restricted to limited scenarios. To address the limitation, we propose the scheme S2. This scheme maintains all the desirable properties and accommodates more general scenarios by splitting the energy and introducing a scalar variable.

Again, we consider the general Wasserstein gradient flow
\begin{equation} \label{wgde3}
    \frac{\partial \rho}{\partial t}= \nabla \cdot\left(\rho \nabla \frac{\delta E}{\delta \rho} \right)
\end{equation}
with a smooth functional $E: \Omega \times \mathbb{R} \rightarrow \mathbb{R}$ and periodic boundary conditions or Neumann boundary conditions. 

Consider the functional \(E\), which can be decomposed into two parts: \(E = E_1 + E_2\). Specifically, 
\[E_1 = E - \int_{\Omega} \rho(\log \rho - 1)\mathrm{~d}\bx,
\text{ and } 
E_2 = \int_{\Omega} \rho(\log \rho - 1)\mathrm{~d}\bx.\]

To employ the SAV method, we assume that \(E_1\) is bounded below,  and there exists a positive constant $C$ such that \(E_1 + C > 0\). Importantly, this assumption is mild in the context of the typical Wasserstein gradient flow. As an illustrative example, we consider \(E(\rho) = \int_{\Omega} \rho^2 dx + 1,\)
which can also be expressed as:
\[E(\rho) = \int_{\Omega} \rho^2 -  \rho(\log \rho - 1) \mathrm{~d} \bx + 1 + \int_{\Omega} \rho(\log \rho - 1) \mathrm{~d} \bx.\]
With 
\(E_1 = \int_{\Omega} \rho^2 -  \rho(\log \rho - 1)\mathrm{~d} \bx + 1,\)
it is evident from the condition \( \rho^2 -  \rho(\log \rho - 1) > 0 \) for \(\rho > 0\) that \(E_1 > 1\).

Building on this foundation, we can introduce a scalar variable $r = \sqrt{E_1+C}$, then the equation \eqref{wgde3} can be expressed as follows:
\begin{equation}
    \frac{\partial \rho}{\partial t}= \nabla \cdot\left(\rho \nabla \left( \frac{r}{ \sqrt{E_1+C}} \frac{\delta E_1}{\delta \rho} + \frac{\delta E_2}{\delta \rho} \right) \right).
\end{equation}
Noting that $\frac{\delta E_2}{\delta \rho}  = \log \rho$, we can deal with $r$ and $\frac{\delta E_2}{\delta \rho}$ implicitly, and obtain:
\begin{align}
         &\frac{\rho^{n+1}- \rho^n}{\delta t} = \nabla \cdot\left(\rho^n \nabla \left( \frac{r^{n+1}}{ \sqrt{E_1^n+C}} \frac{\delta E^n_1}{\delta \rho} + \log \rho^{n+1} \right) \right) \tag{S2-1},\label{sav31}\\
         &\frac{r^{n+1}-r^n}{\delta t} = \frac{1}{2\sqrt{E_1^n+C}}\left(\frac{\delta E^n_1}{\delta \rho},\frac{\rho^{n+1}- \rho^n}{\delta t}\right). \tag{S2-2}\label{sav32}
\end{align}
    \begin{theorem}\label{thm:S2}
    The scheme S2 exhibits the following attributes:
    \begin{itemize}
        \item[1.] Mass conservation:
        $ \int_{\Omega} \rho^n \mathrm{~d} \bx= \int_{\Omega} \rho^{n+1} \mathrm{~d} \bx$;
        \item[2.] Positivity preserving: $\rho^{n+1}>0$ if $\rho^n>0$;
        \item[3.] Unique solvability;
        \item[4.] Energy dissipation law: 
        \begin{equation}\label{s2:energy diss}
     \tilde{E}^{n+1} - \tilde{E}^n \leq -\delta t\int_{\Omega} \rho^n \left|\nabla \left( \frac{r^{n+1}}{ \sqrt{E_1^n+C}} \frac{\delta E^n_1}{\delta \rho} + \log \rho^{n+1} \right)\right|^2 \mathrm{~d}\bx,
\end{equation}
        where the discrete energy $\tilde{E}^n = \int_{\Omega} \rho^{n} \left(\log  \rho^{n} -1 \right) \mathrm{~d} \bx + \left(r^n\right)^2$.
    \end{itemize}
\end{theorem}
\begin{proof}
    Integrating equation \eqref{sav31} over the domain $\Omega$ and applying Neumann or periodic boundary conditions to both $\frac{\delta E_1}{\delta \rho}$ and $\log \rho$ enables us to establish the mass conservation law.

    Furthermore, the inclusion of $\log \rho^{n+1}$ in our formulation guarantees the preservation of positivity for $\rho^{n+1}$.

    Next, we show that (S2) is  uniquely solvable. We first plug the equation \eqref{sav31} into \eqref{sav32}, and obtain
\begin{equation}
     \frac{\rho^{n+1}- \rho^n}{\delta t} = \frac{r^{n} + \frac{1}{2\sqrt{E_1^n+C}}\left(\frac{\delta E^n_1}{\delta \rho},{\rho^{n+1}- \rho^n}\right)}{ \sqrt{E_1^n+C}} \nabla \cdot\left(\rho^n \nabla  \frac{\delta E^n_1}{\delta \rho} \right) + \nabla \cdot\left(\rho^n \nabla \log \rho^{n+1} \right).\label{sav33}
\end{equation}
 We then define $\mathcal{L}_n$ by $\mathcal{L}_n g = u$ where $u$ is the solution of the following elliptic equation
    \begin{equation*}
         - \nabla \cdot \left( \rho^n \nabla  u   \right) = g, \quad \int_{\Omega} u \mathrm{~d}\bx = 0
    \end{equation*}
     with Neumann boundary conditions or periodic boundary conditions.
    Hence, $\mathcal{L}_n$ is self-adjoint and semi-positive definite. 
Similarly, we define a functional 
\begin{equation}
    \begin{aligned}
        F[\rho^{n+1}]&:=\left(\rho^{n+1} (\log \rho^{n+1} -1 ), 1\right) +\frac{1}{2\delta t} \left(\rho^{n+1}-\rho^{n},\mathcal{L}_n\left(\rho^{n+1}-\rho^{n}\right) \right)\\
        &+\left( \frac{r^{n}}{\sqrt{E_1^n+C}}  \frac{\delta E^n_1}{\delta \rho} , \rho^{n+1}\right) + \frac{1}{ {4(E_1^n+C)}}\left(\frac{\delta E^n_1}{\delta \rho},{\rho^{n+1}- \rho^n}\right)^2.
      \end{aligned}
    \end{equation}
    The Euler-Lagrange equation of the above functional under mass conservation is
\begin{equation}\label{sav34}
\begin{aligned}
      &\frac{1}{\delta t} \mathcal{L}_n(\rho^{n+1}- \rho^n) +  \frac{r^{n}}{\sqrt{E_1^n+C}}  \frac{\delta E^n_1}{\delta \rho} + \frac{1}{ 2{E_1^n}}\left(\frac{\delta E^n_1}{\delta \rho},{\rho^{n+1}- \rho^n}\right) \frac{\delta E^n_1}{\delta \rho}  +  \log \rho^{n+1} = \lambda,\\
      &\int_{\Omega} \rho^{n+1} - \rho^{n} \mathrm{~d}\bx = 0,
\end{aligned}
\end{equation}
 where $\lambda$ is the Lagrange multiplier to enforce the mass conservation. It is easy to see from the definition of $\mathcal{L}_n$ that the above is equivalent to \eqref{sav33}.
 
 On the other hand, $F$ is clearly a strictly convex functional on $$\mathcal{A}=\left\{\rho^{n+1} \in H^2(\Omega): \rho^{n+1}>0 \right\}.$$ Hence $F$ admits a unique minimizer, i.e., the equation \eqref{sav34} is  uniquely solvable, which implies that (S2) is uniquely solvable. Furthermore, the minimizer  $\rho^{n+1}$ is positive since the derivative of the term $\rho^{n+1} (\log \rho^{n+1} -1 )$ tends to $-\infty$ at zero.

    To show the energy dissipation, we take the inner product of (S2-1) with $\frac{r^{n+1}}{ \sqrt{E_1^n+C}} \frac{\delta E^n_1}{\delta \rho} + \log \rho^{n+1}$, we obtain
    \begin{equation*}
    \begin{aligned}
        &\left(\frac{r^{n+1}}{ \sqrt{E_1^n+C}} \frac{\delta E^n_1}{\delta \rho} + \log \rho^{n+1},\frac{\rho^{n+1}- \rho^n}{\delta t}\right)\\
         &= -\left( \rho^n \nabla \left( \frac{r^{n+1}}{ \sqrt{E_1^n+C}} \frac{\delta E^n_1}{\delta \rho} + \log \rho^{n+1} \right) , \nabla \left( \frac{r^{n+1}}{ \sqrt{E_1^n+C}} \frac{\delta E^n_1}{\delta \rho} + \log \rho^{n+1} \right)\right).
    \end{aligned}
\end{equation*}
  On the other hand, taking the inner product  of (S2-2) with $2r^{n+1}$, and using a Taylor expansion, we can rewrite the first term in the above to
\begin{equation}
    \begin{aligned}
         &\left(\frac{r^{n+1}}{ \sqrt{E_1^n+C}} \frac{\delta E^n_1}{\delta \rho} + \log \rho^{n+1},\frac{\rho^{n+1}- \rho^n}{\delta t}\right)\\
         &=\frac{2({r^{n+1}-r^n)r^{n+1}}}{\delta t} + \left(\log \rho^{n+1},\frac{\rho^{n+1}- \rho^n}{\delta t}\right)\\
         &=\frac{1}{\delta t}\left\{ \left(r^{n+1}\right)^2 -  \left(r^{n}\right)^2 +  \left(r^{n+1}-r^n\right)^2\right\}\\& \quad +\frac{1}{\delta t}\left\{\left(\log \rho^{n+1} -1 \right) \rho^{n+1} - \left(\log \rho^{n} -1 \right) \rho^{n} + \frac{\left( \rho^{n+1}-\rho^{n}\right)^2}{2\xi}\right\}.
    \end{aligned}
\end{equation}
Combining the above two identities, we have
\begin{equation*}
\begin{aligned}
     &\tilde{E}^{n+1} - \tilde{E}^n\\ & \leq -\delta t\left( \rho^n \nabla \left( \frac{r^{n+1}}{ \sqrt{E_1^n+C}} \frac{\delta E^n_1}{\delta \rho} + \log \rho^{n+1} \right) , \nabla \left( \frac{r^{n+1}}{ \sqrt{E_1^n+C}} \frac{\delta E^n_1}{\delta \rho} + \log \rho^{n+1} \right)\right)\\&\leq 0,
\end{aligned}
\end{equation*}
where $\tilde{E}^n = \int_{\Omega} \rho^{n} \left(\log  \rho^{n} -1 \right) \mathrm{~d} \bx + \left(r^n\right)^2$.
\end{proof}

We can also construct a second-order scheme with BDF and Adams-Bashforth extrapolation as follows:
\begin{align}
         &\frac{3\rho^{n+1}- 4\rho^n+\rho^{n-1}}{2\delta t} = \nabla \cdot\left(\rho^{*,n+\frac{1}{2}} \nabla \left( \frac{r^{n+1}}{ \sqrt{E_1^{*,n+\frac{1}{2}}+C}} \frac{\delta E_1^{*,n+\frac{1}{2}}}{\delta \rho} + \log \rho^{n+1} \right) \right), \label{s2-2nd-1}\\
         &\frac{3r^{n+1}-4r^n+r^{n-1}}{2\delta t} = \frac{1}{2\sqrt{E_1^{*,n+\frac{1}{2}}+C}}\left(\frac{\delta E_1^{*,n+\frac{1}{2}}}{\delta \rho},\frac{3\rho^{n+1}- 4\rho^n+\rho^{n-1}}{2\delta t}\right), \label{s2-2nd-2}
\end{align}
where $E_1^{*,n+\frac{1}{2}} = E_1(\rho^{*,n+\frac{1}{2}})$, $\frac{\delta E_1^{*,n+\frac{1}{2}}}{\delta \rho} = \frac{\delta E_1}{\delta \rho}(\rho^{*,n+\frac{1}{2}})$ and $\rho^{*,n+\frac{1}{2}}$ is defined as the equation \eqref{abe}.
    \begin{theorem}\label{thm:S2-2nd}
     Assume that $\rho^1$ is obtained from the first-order scheme. The above second-order scheme \eqref{s2-2nd-1}-\eqref{s2-2nd-2} exhibits the following attributes:
      \begin{itemize}
        \item[1.] Mass conservation:
        $ \int_{\Omega} \rho^n \mathrm{~d} \bx= \int_{\Omega} \rho^{n+1} \mathrm{~d} \bx$;
        \item[2.] Positivity preserving: $\rho^{n+1}>0$;
        \item[3.] Unique solvability.
    \end{itemize}
\end{theorem}
\begin{proof}
    The proof is similar to that of Theorem \ref{thm:S2}. We can define a slightly different linear operator $\tilde{\mathcal{L}}_n$ as follows:  $\tilde{\mathcal{L}}_n g = u$ is defined by the following elliptic equation:
    \begin{equation*}
         -\nabla \cdot \left( \rho^{*, n+\frac{1}{2}} \nabla  u   \right) = g, \quad \int_{\Omega} u \mathrm{~d}\bx = 0
    \end{equation*}
     with the homogeneous Neumann boundary conditions or periodic boundary conditions.
     We can also define a slightly different convex functional $\tilde{F}$ as follows:
     \begin{equation*}
    \begin{aligned}
        F[\rho^{n+1}]&:=\left(\rho^{n+1} (\log \rho^{n+1} -1 ), 1\right) +\left( \frac{4r^{n}-r^{n-1}}{3\sqrt{E_1^{*,n+\frac{1}{2}}+C}} \frac{\delta E_1^{*,n+\frac{1}{2}}}{\delta \rho} , \rho^{n+1}\right)\\&+\frac{1}{12\delta t} \left(3\rho^{n+1}-4\rho^{n}+\rho^{n-1},\tilde{\mathcal{L}}_n\left(3\rho^{n+1}-4\rho^{n}+\rho^{n-1}\right) \right)\\
        & + \frac{1}{ {36(E_1^{*,n+\frac{1}{2}}+C)}}\left(\frac{\delta E_1^{*,n+\frac{1}{2}}}{\delta \rho},3\rho^{n+1}-4\rho^{n}+\rho^{n-1}\right)^2.
      \end{aligned}
    \end{equation*}
    Then it can be shown that the solution of \eqref{s2-2nd-1}-\eqref{s2-2nd-2} is the unique minimizer of the above convex functional.
\end{proof}

\subsection{Application to Onsager gradient flows}
The Onsager gradient flow, as presented in \cite{doi2011onsager}, is  as follows:
\begin{equation}\label{Onsage gf}
    \frac{\partial \rho}{\partial t}=\nabla \cdot\left(V_1(\rho) \nabla \frac{\delta  {E}(\rho)}{\delta \rho}\right)-V_2(\rho) \frac{\delta  {E}(\rho)}{\delta \rho},
\end{equation}
where ${E}(\rho)$ is an energy functional, and $V_1, V_2: \mathbb{R} \rightarrow \mathbb{R}_{+}$are positive mobility functions.
By  decomposing the energy functional \( {E}(\rho) \) into two components: \( {E}(\rho) = {E_1}(\rho) + {E_2}(\rho) \), where \( E_1 \) is bounded from below with a constant, \( E_2 \) is convex with respect to \( \rho \), and introducing a scalar variable, \( r \), defined as \( r = \sqrt{{E_1} + C} \), where $C$ is a constant such that ${{E_1} + C}>0$, the equation \eqref{Onsage gf} can be written in the following manner:
\begin{equation*}
\begin{aligned}
    \frac{\partial \rho}{\partial t}&=\nabla \cdot\left(V_1 \nabla \left(\frac{r} {\sqrt{{E_1+C}}}\frac{\delta  {E_1}}{\delta \rho} + \frac{\delta  {E_2}}{\delta \rho}\right)\right)-V_2 \left(\frac{r} {\sqrt{{E_1+C}}}\frac{\delta  {E_1}}{\delta \rho} + \frac{\delta  {E_2}}{\delta \rho}\right).
\end{aligned}
\end{equation*}
Then, following the construction strategy of S1 and denoting $V^n_1 = V_1(\rho^n)$ and $V^n_2 = V_2(\rho^n)$, we construct the following scheme:
\begin{align}
         \frac{\rho^{n+1}- \rho^n}{\delta t} &= \nabla \cdot\left(V^n_1 \nabla \left( \frac{r^{n+1}}{ \sqrt{E_1^n+C}} \frac{\delta E^n_1}{\delta \rho} +  \frac{\delta E^{n+1}_2}{\delta \rho} \right) \right)\label{Onsage gf s31} \\&+ V^n_2 \left( \frac{r^{n+1}}{ \sqrt{E_1^n+C}} \frac{\delta E^n_1}{\delta \rho} +  \frac{\delta E^{n+1}_2}{\delta \rho} \right), \nonumber\\
         \frac{r^{n+1}-r^n}{\delta t} &= \frac{1}{2\sqrt{E_1^n+C}}\left(\frac{\delta E^n_1}{\delta \rho},\frac{\rho^{n+1}- \rho^n}{\delta t}\right) \label{Onsage gf s32}.
\end{align}
Taking the inner products of  \eqref{Onsage gf s31}  with $\frac{r^{n+1}}{ \sqrt{E_1^n+C}} \frac{\delta E^n_1}{\delta \rho} + \frac{\delta E^{n+1}_2}{\delta \rho}$, of \eqref{Onsage gf s32}  with  $2r^{n+1}$, we obtain:
\begin{equation*}
    \begin{aligned}
         & \left(r^{n+1}\right)^2 -  \left(r^{n}\right)^2 +  \left(r^{n+1}-r^n\right)^2 + E^{n+1}_2- E^{n}_2 \\
         &\leq -\delta t\left(V_1^n \nabla \left( \frac{r^{n+1}}{ \sqrt{E_1^n+C}} \frac{\delta E^n_1}{\delta \rho} +  \frac{\delta E^{n+1}_2}{\delta \rho} \right) , \nabla \left( \frac{r^{n+1}}{ \sqrt{E_1^n+C}} \frac{\delta E^n_1}{\delta \rho} +  \frac{\delta E^{n+1}_2}{\delta \rho}\right)\right)\\
         &-\delta t\left( V_2^n \left( \frac{r^{n+1}}{ \sqrt{E_1^n+C}} \frac{\delta E^n_1}{\delta \rho} +  \frac{\delta E^{n+1}_2}{\delta \rho} \right) , \left( \frac{r^{n+1}}{ \sqrt{E_1^n+C}} \frac{\delta E^n_1}{\delta \rho} +  \frac{\delta E^{n+1}_2}{\delta \rho}\right)\right)\leq 0.
    \end{aligned}
\end{equation*}
Through appropriate choice of the form of \( E_2 \), we can impose specific attributes for the scheme \eqref{Onsage gf s31}-\eqref{Onsage gf s32}. As an illustrative example, by selecting \( E_2(\rho) = \int_{\Omega} \rho(\log \rho - 1) \mathrm{~d} \boldsymbol{x}\) , we can ensure that \( \rho \) remains  positive.
\begin{theorem}
    With \( E_2(\rho) = \int_{\Omega} \rho(\log \rho - 1) \mathrm{~d} \boldsymbol{x}\), the scheme \eqref{Onsage gf s31}-\eqref{Onsage gf s32}, exhibits the following attributes:
    \begin{itemize}
        \item[1.] Positivity preserving: $\rho^{n+1}>0$;
        \item[2.] Unique solvability;
        \item[3.] Energy dissipation law: 
        \begin{equation}
            \tilde{E}^{n+1} - \tilde{E}^n \leq 0,
        \end{equation}
        where the discrete energy is defined as $\tilde{E}^n = \int_{\Omega} \rho^{n} \left(\log  \rho^{n} -1 \right) \mathrm{~d} \bx + \left(r^n\right)^2$.
    \end{itemize}
\end{theorem}
\begin{proof}
    The proof is almost the same as the Theorem \ref{thm:S2} excluding the mass conservation.
\end{proof}
\section{Finite difference discretization in space}
In this section, we turn our attention to the construction of finite difference schemes that can maintain the properties of the time discretizations discussed in the previous section. Since integration by parts is essential in mass conservation and energy dissipation in the spatial continuous case,  
a critical aspect in full discretization is to make sure that the implementation of boundary conditions satisfies suitable  summation by parts formulae. This task is relatively straightforward in rectangular domains.  Since the construction methodologies for S1 are essentially the same as those presented in \cite{ShenXu2021}. Our subsequent discussion primarily focuses on the 2D implementation of S2. 

Let the domain $[0, L]^2$ be discretized into $M^2$ points, designated as \\
$x_{j, k} = \left(\left(j-\frac{1}{2}\right) \delta x, \left(k-\frac{1}{2}\right) \delta x\right)$ for $j, k = 1, \ldots, M$, with $\delta x = L / M$.
Then, a fully discrete version of the  scheme (S2-1)-(S2-2) is 
\begin{equation}\label{sav31fd}
    \begin{aligned}
    \frac{\rho_{j, k}^{n+1}-\rho_{j, k}^{n}}{\delta t}&=\frac{1}{\delta x^{2}}\left[\frac{\rho_{j+1, k}^{n}+\rho_{j, k}^{n}}{2}\left(\xi^{n+1}\varphi_{j+1, k}^{n} + \mu_{j+1, k}^{n+1}-\xi^{n+1}\varphi_{j, k}^{n}-\mu_{j, k}^{n+1}\right)\right.\\
&-\frac{\rho_{j, k}^{n}+\rho_{j-1, k}^{n}}{2}\left(\xi^{n+1}\varphi_{j, k}^{n} + \mu_{j, k}^{n+1}-\xi^{n+1}\varphi_{j-1, k}^{n}-\mu_{j-1, k}^{n+1}\right), \\
&+\frac{\rho_{j, k+1}^{n}+\rho_{j, k}^{n}}{2}\left(\xi^{n+1}\varphi_{j, k+1}^{n} + \mu_{j, k+1}^{n+1}-\xi^{n+1}\varphi_{j, k}^{n}-\mu_{j, k}^{n+1}\right) \\
&\left.-\frac{\rho_{j, k}^{n}+\rho_{j, k-1}^{n}}{2}\left(\xi^{n+1}\varphi_{j, k}^{n} + \mu_{j, k}^{n+1}-\xi^{n+1}\varphi_{j, k-1}^{n}-\mu_{j, k-1}^{n+1}\right)\right],
\end{aligned}
\end{equation}
\begin{equation}\label{sav32fd}
    {r^{n+1}-r^{n}}=\frac{\delta x^2}{2\sqrt{E_1^n+C}}\sum_{j,k=1}^{M} \varphi_{j, k}^{n}\left(\rho_{j, k}^{n+1}-\rho_{j, k}^{n}\right),
\end{equation}
where $\xi^{n+1} = \frac{r^{n+1}}{ \sqrt{E_1^n+C}}, \varphi_{j, k}^{n}=  (\frac{\delta E_1}{\delta \rho})_{j, k}^{n} \text{ and } \mu_{j, k}^{n} = (\log \rho)_{j, k}^{n} \text{ for } 1 \leq j, k \leq M$. To illustrate the handling of boundary conditions, we consider, as an example,  the Neumann boundary conditions. For achieving summation by parts, it is necessary to impose boundary terms as follows:
\begin{equation}\label{bdc}
\begin{aligned}
& \frac{\varphi_{0, k}^{n+1}-\varphi_{1, k}^{n+1}}{\delta x}=0, \quad \frac{\varphi_{M+1, k}^{n+1}-\varphi_{M, k}^{n+1}}{\delta x}=0,\\
& \frac{\mu_{0, k}^{n+1}-\mu_{1, k}^{n+1}}{\delta x}=0, \quad \frac{\mu_{M+1, k}^{n+1}-\mu_{M, k}^{n+1}}{\delta x}=0.
\end{aligned}
\end{equation}

By representing $\tilde{\rho}^n$ as a vector composed of the elements $\rho_{j, k}^{n}$, $\tilde{\varphi}^n$ as a vector consisting of $\varphi_{j, k}^{n}$, and $\tilde{\mu}^n$ as a vector formed from $\mu_{j, k}^{n+1}$ for $1 \leq j, k \leq M$ arranged in lexicographical order, we can rewrite the above into a matrix form as follows:
\begin{equation}\label{s3mat1}
\begin{aligned}
    \frac{\delta x^2}{\delta t} \left(\tilde{\rho}^{n+1} - \tilde{\rho}^n\right) = -A^n \left(\frac{r^n}{{\sqrt{E_1^n + C}}} \tilde{\varphi}^n + \frac{\delta x^2}{2{(E_1^n + C)}} (\tilde{\varphi}^n)^T \left(\tilde{\rho}^{n+1} - \tilde{\rho}^n\right)\tilde{\varphi}^n + \tilde{\mu}^{n+1} \right).
\end{aligned}
\end{equation}
Here, $A^n$ is a sparse matrix with nonzero elements adjacent to its diagonal. Two indices $(j,k)$ and $(j',k')$ are considered adjacent if $\left|j-j'\right| + \left|k-k'\right| = 1$. The diagonal entry of $A^n$ is:
\begin{equation}
\left[2\rho_{j, k}^n + \frac{1}{2} \left(\rho_{j+1, k}^n + \rho_{j-1, k}^n + \rho_{j, k+1}^n + \rho_{j, k-1}^n\right)\right].
\end{equation}
$A^n$ is symmetric and diagonally dominant, with positive diagonal entries, ensuring it is positive semi-definite with non-negative eigenvalues. The eigen-decomposition of $A^n$ is $A^n = T^t \Lambda T$ where $\Lambda = \operatorname{diag}\left(0, \mu_2, \ldots, \mu_{M^2}\right)$ and $\mu_j > 0$ for $j = 2, \ldots, M^2$. The pseudo-inverse of $A^n$, denoted $(A^n)^*$, is defined as $(A^n)^*= T^t \operatorname{diag}\left(0, \mu_2^{-1}, \ldots, \mu_{M^2}^{-1}\right) T$. For the term $(\tilde{\varphi}^n)^T\left(\tilde{\rho}^{n+1} - \tilde{\rho}^n\right)\tilde{\varphi}^n$, we have:
\begin{equation*}
(\tilde{\varphi}^n)^T\left(\tilde{\rho}^{n+1} - \tilde{\rho}^n\right)\tilde{\varphi}^n = \tilde{\varphi}^n (\tilde{\varphi}^n)^T\left(\tilde{\rho}^{n+1} - \tilde{\rho}^n\right).
\end{equation*}
Defining $\Psi^n = \tilde{\varphi}^n (\tilde{\varphi}^n)^T$, a positive semi-definite matrix, and multiplying equation \eqref{s3mat1} by the pseudo-inverse $(A^n)^*$, we obtain:
\begin{equation}\label{s3mat2}
\frac{\delta x^2}{\delta t} (A^n)^*\left(\tilde{\rho}^{n+1} - \tilde{\rho}^n\right) + a_1^n \tilde{\varphi}^n + a_2^n \Psi^n \left(\tilde{\rho}^{n+1} - \tilde{\rho}^n\right) + \tilde{\mu}^{n+1} = 0,
\end{equation}
where $a_1^n = \frac{r^n}{{\sqrt{E_1^n + C}}}$ and $a_2^n = \frac{\delta x^2}{2(E_1^n + C)}$.

For the above scheme, we have
\begin{theorem}\label{thms3}
    The finite difference scheme \eqref{sav31fd}-\eqref{sav32fd} enjoys the following properties:
        \begin{itemize}
        \item[(i)] Mass conservation:$$
\delta x^2 \sum_{j, k=1}^M \rho_{j, k}^{n+1} = \delta x^2 \sum_{j, k=1}^M \rho_{j, k}^n, 1 \leq i \leq N;
$$
        \item[(ii)] Positivity preserving: if \(\rho_{j, k}^n > 0\) for \((j, k)\), we have \(\rho_{j, k}^{n+1} > 0\) for all \((j, k)\);
        \item[(iii)] Unique solvability: the scheme \eqref{sav31fd} possesses a unique solution $\rho_{j, k}^{n+1}$;
        \item[(iv)] Energy dissipation law: 
$$
\begin{aligned}
&\tilde{E}^{n+1} - \tilde{E}^n\\ &\leq  -\frac{\delta t }{\delta x^2} \sum_{\substack{1 \leq j \leq M-1 \\
1 \leq k \leq M}} \frac{\rho_{j+1, k}^n + \rho_{j, k}^n}{2}\left(\xi^{n+1}\varphi_{j+1, k}^{n} + \mu_{j+1, k}^{n+1}-\xi^{n+1}\varphi_{j, k}^{n}-\mu_{j, k}^{n+1}\right)^2 \\
& +\sum_{\substack{1 \leq j \leq M \\
1 \leq k \leq M-1}} \frac{\rho_{j, k+1}^n + \rho_{j, k}^n}{2}\left(\xi^{n+1}\varphi_{j, k+1}^{n} + \mu_{j, k+1}^{n+1}-\xi^{n+1}\varphi_{j, k}^{n}-\mu_{j, k}^{n+1}\right)^2,
\end{aligned}
$$
where \begin{equation*}
    \tilde{E}^n = \sum_{j,k=1}^{M} \rho_{j, k}^n\left(\log \rho_{j, k}^n -1 \right) + (r^n)^2.
\end{equation*}
\end{itemize}
\end{theorem}
\begin{proof}
    The mass conservation is obtained by taking the sum over $1 \leq j, k \leq M$ on equation \eqref{sav31fd} and using the boundary conditions of $\left(\varphi\right)_{j, k}^{n}$ and $\left(\mu\right)_{j, k}^{n+1}$ in \eqref{bdc}.

To prove the unique solvability and positivity preservation, we define  
\begin{equation}\label{convexF}
\begin{aligned}
     F[\tilde{\rho}^{n+1}] &= \frac{\delta x^2}{2\delta t} \left(\tilde{\rho}^{n+1} - \tilde{\rho}^n\right)^T(A^n)^*\left(\tilde{\rho}^{n+1} - \tilde{\rho}^n\right) + a_1^n \left(\tilde{\varphi}^n\right)^T \tilde{\rho}^{n+1} \\&+ \frac{a_2^n}{2} \left(\tilde{\rho}^{n+1} - \tilde{\rho}^n\right)^T\Psi^n \left(\tilde{\rho}^{n+1} - \tilde{\rho}^n\right) + \left(\tilde{\rho}^{n+1}\right)^T\left(\log \tilde{\rho}^{n+1} -1\right).
\end{aligned}
\end{equation}
It is evident that $F[\tilde{\rho}^{n+1}]$ is strictly convex with respect to $\tilde{\rho}^{n+1}$, and equation \eqref{s3mat2} represents its Euler-Lagrange equation. Consequently, the solution to \eqref{s3mat2} is the unique minimizer of $F[\tilde{\rho}^{n+1}]$ with $\tilde{\rho}^{n+1} > 0$. As any element of $\tilde{\rho}^{n+1}$ approaches zero, the gradient of $F[\tilde{\rho}^{n+1}]$ tends towards negative infinity, implying that the function's value decreases when incrementing elements near 0 in $\tilde{\rho}^{n+1}$. Therefore, a minimizer with any zero element in $\tilde{\rho}^{n+1}$ is not attainable.

The energy dissipation law can be obtained by multiplying the equation \eqref{sav31fd} with $\xi^{n+1}\varphi_{j, k}^{n}+\mu_{j, k}^{n+1}$ and multiplying the equation \eqref{sav32fd} with $2r^{n+1}$. The summation by parts is used with the employed boundary conditions of $\varphi_{j, k}^{n}$ and $\mu_{j, k}^{n+1}$.
\end{proof}

\section{Numerical experiments}
In this section, we present various numerical experiments  to validate the theoretical results discussed previously. It is important to note that  our schemes require solving a nonlinear system at every time step. To solve these nonlinear equations, we employ the damped Newton’s iteration method \cite{inbook}. In scenarios where the density must remain positive, owing to the presence of logarithmic functions in the equations, we implement a small corrective perturbation. This is achieved by setting $\rho = \max(\rho, \epsilon)$, thereby ensuring the density does not fall below a minimal threshold, denoted by $\epsilon$. We set $\epsilon=10^{-6}$ in the following experiments.
\subsection{Accuracy test}
First, we test the accuracy in time of the scheme S1 and S2 by solving a heat equation with Neumann boundary conditions:
\begin{equation}\label{heat equation}
\begin{aligned}
&\frac{\partial \rho}{\partial t} = \frac{1}{50}\rho_{xx}=\frac{1}{50} \partial_x \rho\partial_x \log(\rho), \quad x \in [0,1], t>0,\\
&\rho_x(0, t) = \rho_x(1,t) = 0,
\end{aligned}
\end{equation}
which can be expressed as a Wasserstein gradient flow with  $E(\rho) =\int_{\Omega}  \frac{1}{50}\rho(\log \rho - 1) \mathrm{~d} \boldsymbol{x}$.
The exact solution is taken as $\rho(x,t) = e^{-\pi^2t/50 }\cos(\pi x) + 1.1$. To assess the accuracy of our model, we employ both the \(L_{\infty}\) error and the \(L_2\) error. These are defined by:
\begin{align*}
e_{\infty}^N &:= \max _i \left| \rho_i^N - \rho\left(x_i, T\right) \right|, \\
e_2^N &:= \left( \delta x \sum_{i=1}^{I-1} \left( \rho_i^N - \rho\left(x_i, T\right) \right)^2 \right)^{\frac{1}{2}},
\end{align*}
where $\rho_i^N$ is the value of $\rho^N$ at $x_i$.

To determine the convergent rate of different schemes, we utilize a fine mesh with a finite difference method, specifically setting \(N = 50000\). The time steps chosen are $\delta t = 0.1, 0.05, 0.025, 0.0125$, with a final time of $T=1$. For the scheme S2, we decompose the energy as
$$E(\rho) =\int_{\Omega}  \frac{1}{50}\rho(\log \rho - 1) \mathrm{~d} \boldsymbol{x}= \int_{\Omega}  \frac{1}{100}\rho(\log \rho - 1) \mathrm{~d} \boldsymbol{x}+ \int_{\Omega}  \frac{1}{100}\rho(\log \rho - 1) \mathrm{~d} \boldsymbol{x}.$$ 
It is observed that the schemes both S1 and S2 consistently demonstrate first-order convergent rates in time.
\begin{table}[tbh]
    \centering
    $\begin{array}{||c|ccccc||}
\hline
\delta x = 1/50000 & & e_{\infty}^N & e_2^N & \text{order of }e_{\infty}^N & \text{order of }e_2^N\\
\hline  \delta t = 0.1 &  & 8.1540e-03 &  3.0473e-03 & - & - \\ 
\hline \delta t = 0.05  &  &   4.1101e-03 & 1.5456e-03 & 0.9883 & 0.9794 \\ 
\hline \delta t = 0.025  &  &  2.0578e-03 & 7.7746e-04 & 0.9981 & 0.9913 \\ 
\hline \delta t = 0.0125  &  & 1.0244e-03 & 3.8890e-04 & 1.0063 & 0.9994 \\
\hline
\end{array}$
    \caption{Heat equation, the first-order convergent rate of S1 in time.}
    \label{tab:Heat equation l2s1}
\end{table}

\begin{table}[tbh]
    \centering
    $\begin{array}{||c|ccccc||}
\hline
\delta x = 1/50000 & & e_{\infty}^N & e_2^N & \text{order of }e_{\infty}^N & \text{order of }e_2^N\\
\hline  \delta t = 0.1 &  & 3.2798e-03 &  1.3060e-03 & - & - \\ 
\hline \delta t = 0.05  &  &  1.6497e-03  & 6.5815e-04 & 0.9794 & 0.9794 \\ 
\hline \delta t = 0.025   &  &  8.2241e-04 & 3.2926e-04 & 0.9913 & 0.9913 \\ 
\hline \delta t = 0.0125  &  & 4.0556e-04 & 1.6352e-04 & 0.9994 & 0.9994 \\
\hline
\end{array}$
    \caption{Heat equation, the first-order convergent rate of S2 in time.}
    \label{tab:Heat equation s3}
\end{table}
\subsection{Barenblatt Solution}
The Barenblatt solution, a fundamental benchmark for the porous medium equation (PME), is widely used to test the accuracy and efficiency of numerical methods. The PME, denoted by \( \frac{\partial \rho}{\partial t} = \Delta \rho^m \), integrates energy as \( E = \int_{\Omega} \frac{1}{m-1} \rho^m(\bx) \, \mathrm{~d}\bx \) and takes the explicit form:
\[
B_{m, d}(\bx, t) = (t+1)^{-\alpha} \left[ 1 - \frac{\alpha(m-1)}{2md} \frac{\|\bx\|^2}{(t+1)^{2\alpha/d}} \right]_+^{1/(m-1)},
\]
where \( (s)_+ = \max(s, 0) \) and \( \alpha = \frac{d}{d(m-1)+2} \). In this context, \( d=2 \) represents the spatial dimensions of the problem. To illustrate the precision and energy dissipation efficiency of our proposed numerical schemes, we evaluate the solution from \( t_0 = 0 \) to \( T = 1 \), with a time step \( \delta t=0.001 \) and spatial step \( \delta x = 0.25 \), over the domain \( \Omega = (-10, 10)^2 \). We utilize schemes S1 and S2 for \( m = 3 \). For the scheme S1, the energy is defined as:
\[
E_{S1} = \int_{\Omega} \rho(\bx)\left(\log \rho(\bx) - 1\right) \mathrm{d}\bx.
\]
For the scheme S2, the original energy is defined as:
\[
E_{S2} = \int_{\Omega} \frac{1}{m-1} \rho^m(\bx) \mathrm{~d}\bx,
\]
and \(E_{S2}\) is decomposed into \(E_{S2} = E_{1,S2} + \int_{\Omega} \rho(\bx)\left(\log \rho(\bx) - 1\right) \mathrm{~d}\bx\).
The modified energy for S2 is then defined as:
\[
\tilde{E}_{S2} = r^2 + \int_{\Omega} \rho(\bx)\left(\log \rho(\bx) - 1\right) \mathrm{~d}\bx.
\]
It is demonstrated that the quantity \(E_{1,S2}\), defined by
\[
E_{1,S2} = \int_{\Omega} \frac{\rho^m(\bx)}{m-1} - \rho(\bx)\left(\log \rho(\bx) - 1\right) \mathrm{~d}\bx,
\]
is bounded from below when \(m \geq 2\). In this context, setting \(C = 0\) for S2 is sufficient to ensure \(E_{1,S2}\) remains positive.

\begin{figure}[tbh]
\centering
     \begin{subfigure}{0.48\textwidth}
         \centering
         \includegraphics[width=\textwidth]{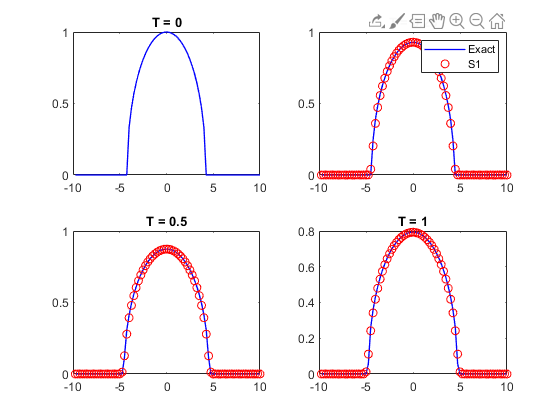}
         \caption{S1}
         \label{fig:pmem3s1rho}
     \end{subfigure}
     \begin{subfigure}{0.48\textwidth}
         \centering
         \includegraphics[width=\textwidth]{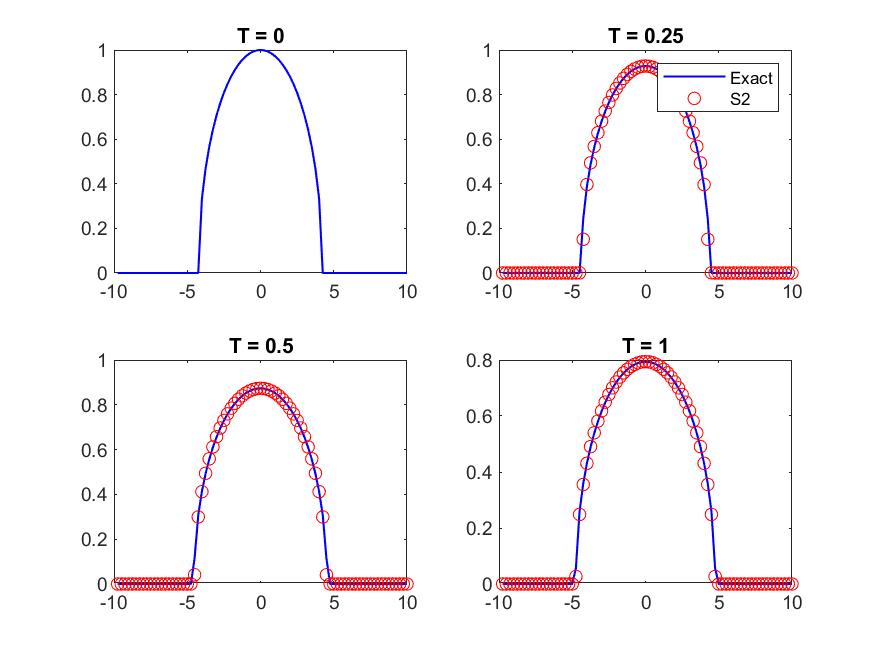}
         \caption{S2}
         \label{fig:pmem3s3rho}
     \end{subfigure}
        \begin{subfigure}{0.48\textwidth}
         \centering
         \includegraphics[width=\textwidth]{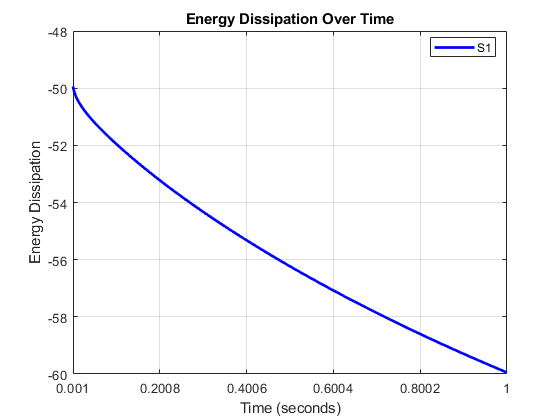}
         \caption{Energy dissipation of S1}
         \label{fig:pmem3s1e}
     \end{subfigure}
     \begin{subfigure}{0.48\textwidth}
         \centering
         \includegraphics[width=\textwidth]{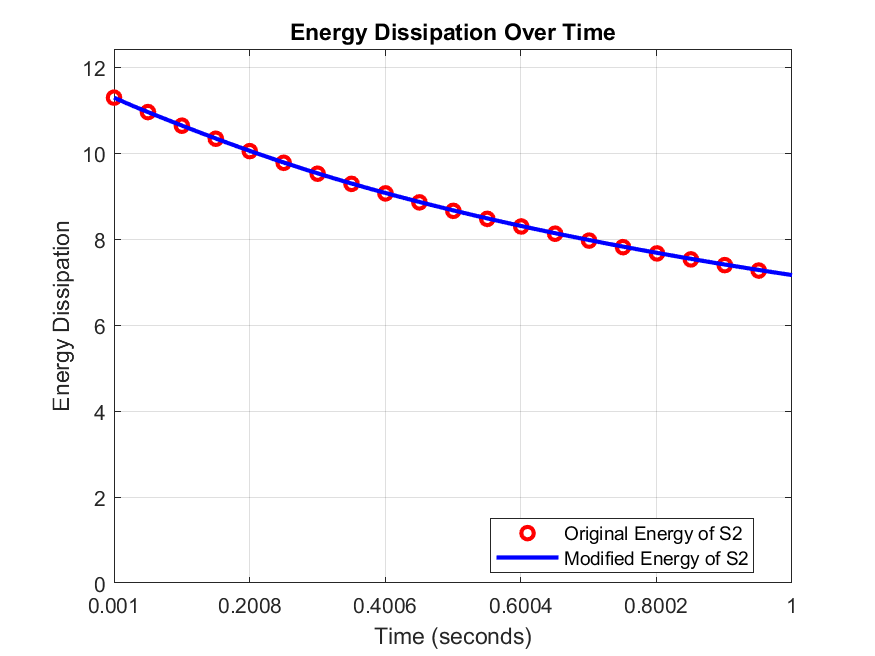}
         \caption{Energy dissipation of S2}
         \label{fig:pmem3s3e}
     \end{subfigure}
     \caption{Visualizations of the PME evolution with S1 and S2 at a cross section \(y=0\) and energy dissipation comparison for \(m=3\).}
     \label{fig:pmes3m3}
     
\end{figure}

Figures \ref{fig:pmem3s1rho} to \ref{fig:pmem3s3rho} display the solution's evolution at the cross section \( y=0 \) for S1 and S2 respectively. These figures confirm that our schemes can closely approximate the exact solution, demonstrating their accuracy. Figures \ref{fig:pmem3s1e} to \ref{fig:pmem3s3e} illustrate the energy dissipation for S1 and S2 respectively, with Figure \ref{fig:pmem3s3e} also comparing the original and modified energy, which are observed to be consistent.

\begin{figure}[tbh]
     \centering
         \includegraphics[width=0.7\textwidth]{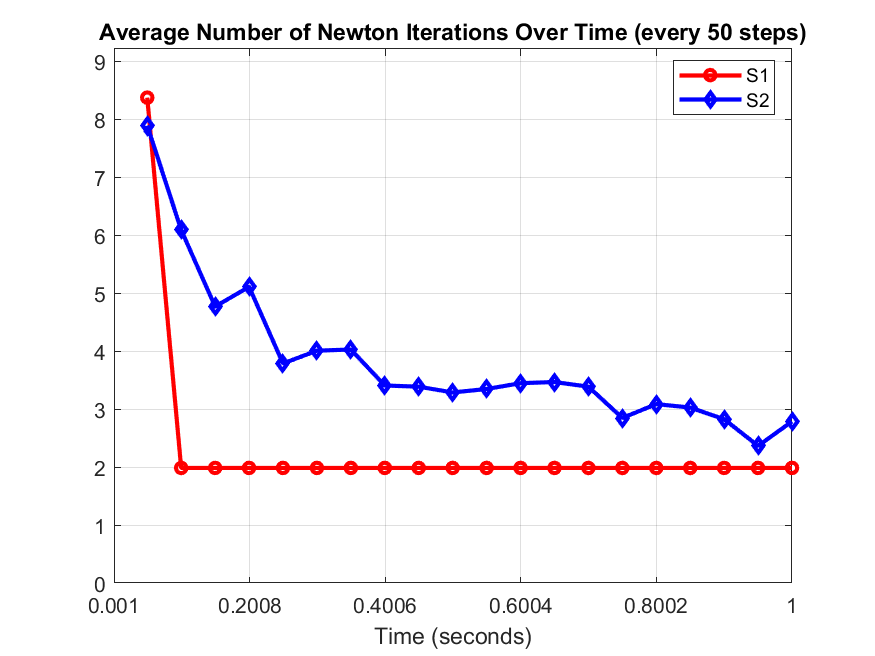}
         \caption{Comparative analysis of the average number of Newton iterations required for schemes S1 and S2 every 50 computational steps over the simulation time from $t_0=0$ to $T=1$, demonstrating the relative computational efficiency of each scheme.}
         \label{fig:pmeavgnewton}
\end{figure}

To evaluate the computational efficiency of our schemes, we compared the average number of Newton iterations required every 50 steps (Figure \ref{fig:pmeavgnewton}). Both S1 and S2 achieve convergence with several Newton iterations, where only the initial few steps require slightly more iterations. Thanks to its simpler structure, the scheme S1 generally needs fewer iterations, making it the preferred option when feasible.

\subsection{Fokker-Planck Equation}
To test the proposed schemes in the context of potential-influenced equations, we turn our attention to the linear Fokker-Planck equation presented as:
\begin{equation*}
\frac{\partial \rho}{\partial t} = \Delta \rho + \nabla \cdot (\rho \nabla V)= \nabla \cdot \rho \nabla (\log \rho+V),
\end{equation*}
which is a Wasserstein gradient flow with the energy functional \(E(\rho)\) as:
\begin{equation*}
E(\rho) = \int_{\Omega} \rho (\log (\rho) - 1) + \rho V \mathrm{~d}\bx.
\end{equation*}
In this example, we take $V(x,y) = \frac{x^2 + y^2}{2}$ which 
will guide the system towards a unique and globally stable equilibrium, irrespective of the initial conditions. The equilibrium state, also known as the heat kernel, is analytically defined as
\begin{equation}\label{heat kernel}
\rho(x,y,t) = \frac{1}{4\pi t} \exp \left( -\frac{x^2 + y^2}{4t} \right),
\end{equation}
evaluated at \(t = \frac{1}{2}\), i.e., $\rho_{\infty} = \rho(x,y,1/2)$.

\begin{figure}[tbh]
 \begin{subfigure}{0.32\textwidth}
         \centering
         \includegraphics[width=\textwidth]{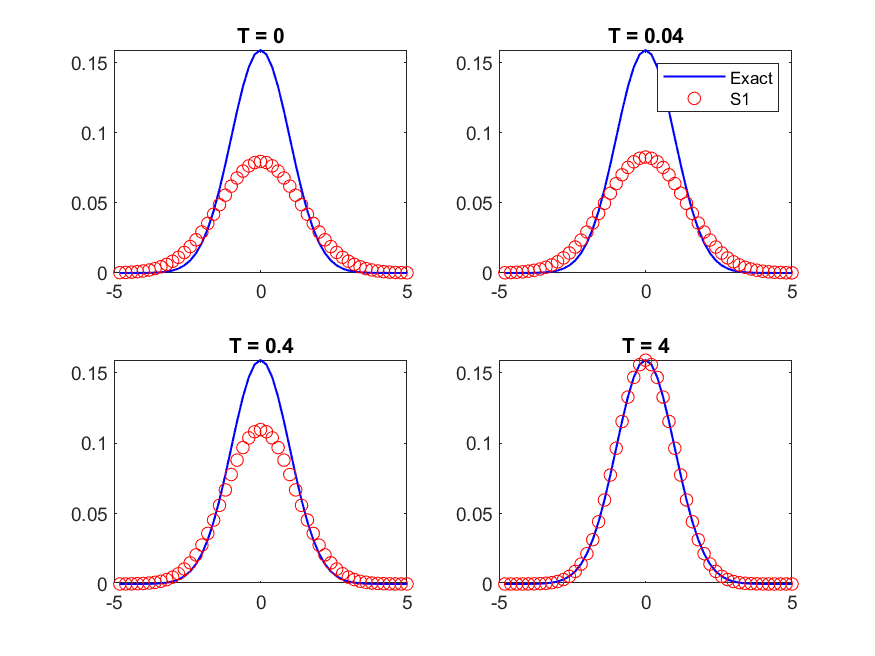}
         \caption{Cross section of S1}
         \label{fkcrosss1}
     \end{subfigure}
    \begin{subfigure}{0.32\textwidth}
         \centering
    \includegraphics[width=\textwidth]{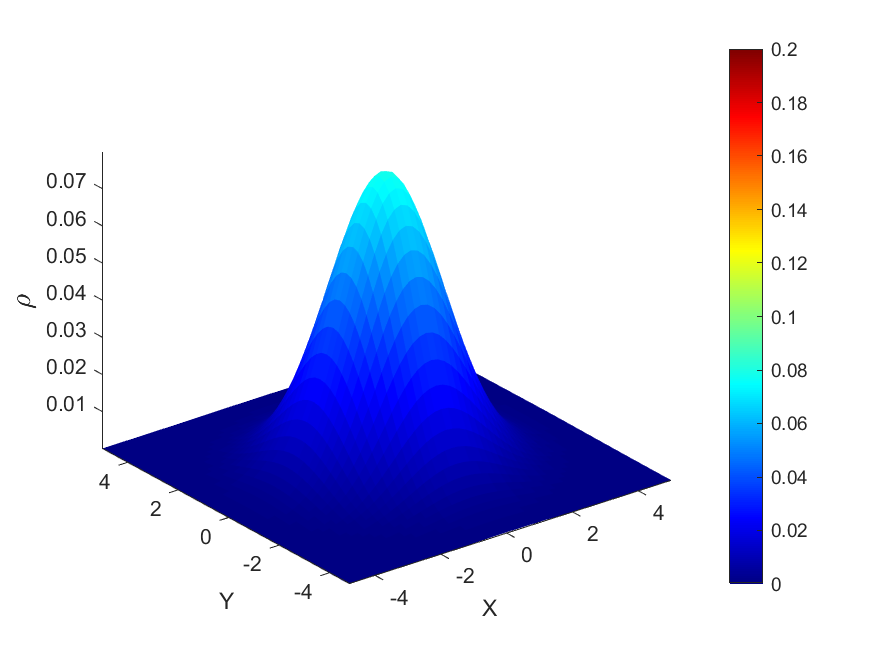}
         \caption{Initial state}
     \end{subfigure}
     \centering
    \begin{subfigure}{0.32\textwidth}
         \centering
         \includegraphics[width=\textwidth]{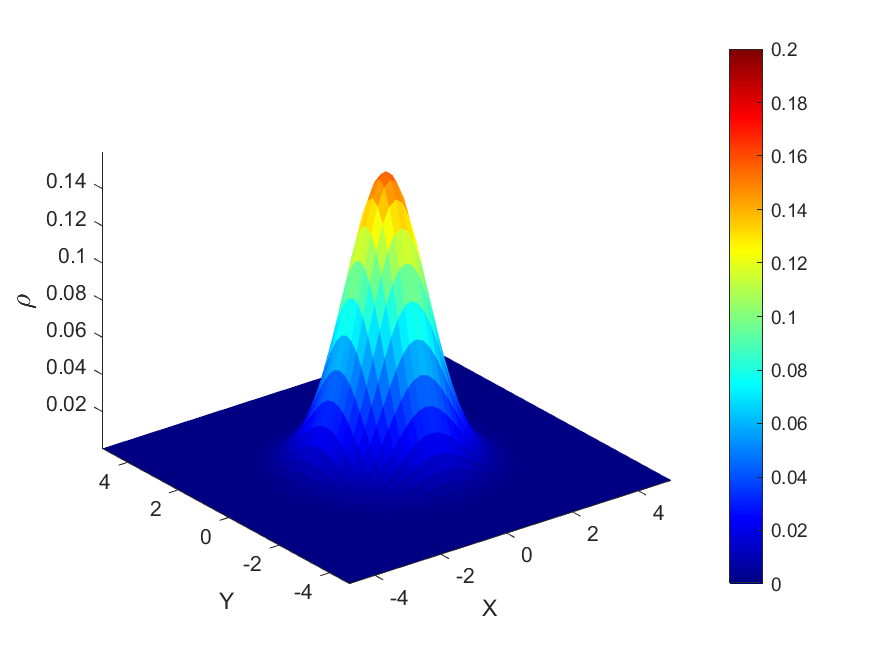}
         \caption{Steady state of S1}
     \end{subfigure}
     \begin{subfigure}{0.32\textwidth}
         \centering
\includegraphics[width=\textwidth]{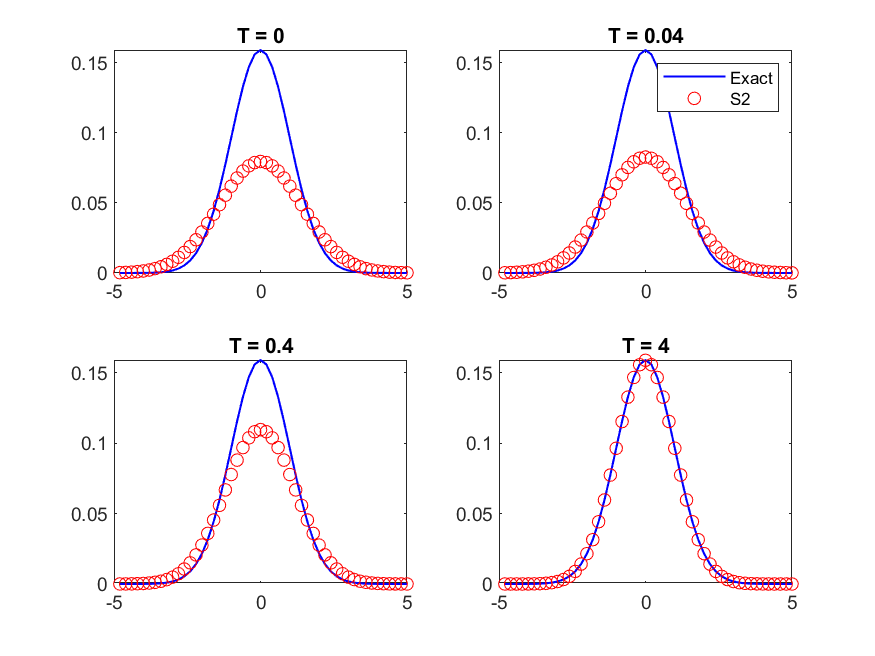}
          \caption{Cross section of S2}
          \label{fkcrosss2}
     \end{subfigure}
          \begin{subfigure}{0.32\textwidth}
         \centering
         \includegraphics[width=\textwidth]{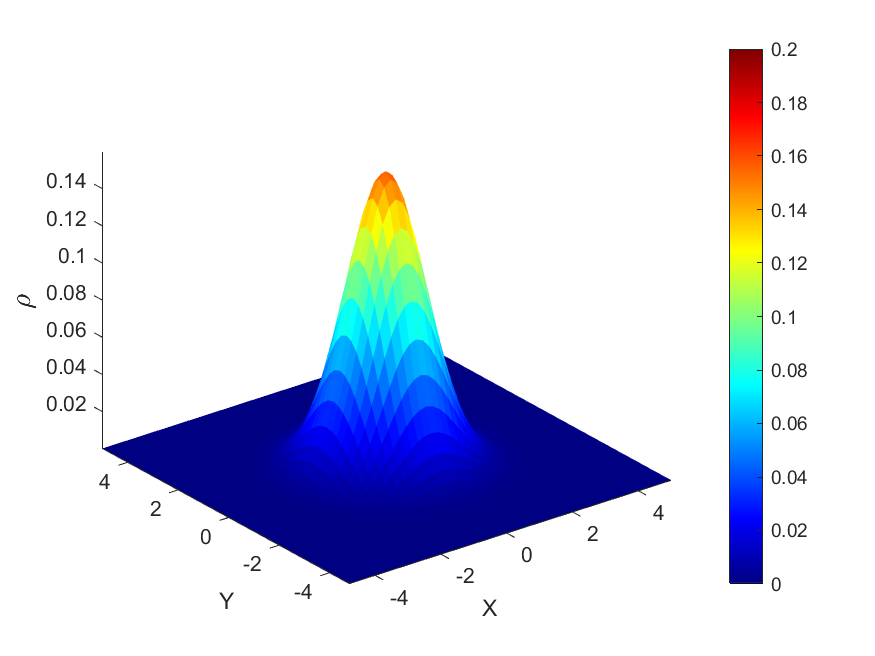}
         \caption{Reference steady state}
     \end{subfigure}
     \begin{subfigure}{0.32\textwidth}
         \centering
         \includegraphics[width=\textwidth]{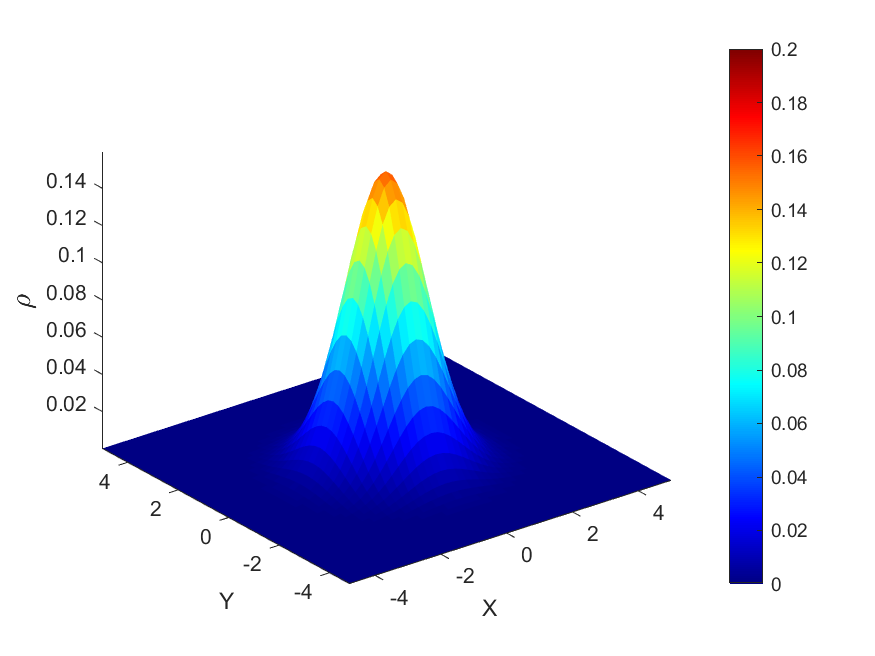}
         \caption{Steady state of S2}
     \end{subfigure}
    \caption{Solution profiles for the Fokker-Planck equation. The simulation begins with the initial condition \(\rho(x,y,1)\) as described in the equation \eqref{heat kernel}. The reference steady state is given by \(\rho_{\infty}\). Additionally, cross-sectional views of the solution at \(y=0\) are presented in Figure \ref{fkcrosss1} and \ref{fkcrosss2}.}
     \label{pmevrand2dinfty2}
\end{figure}

\begin{figure}[tbh]
     \centering
         \includegraphics[width=0.7\textwidth]{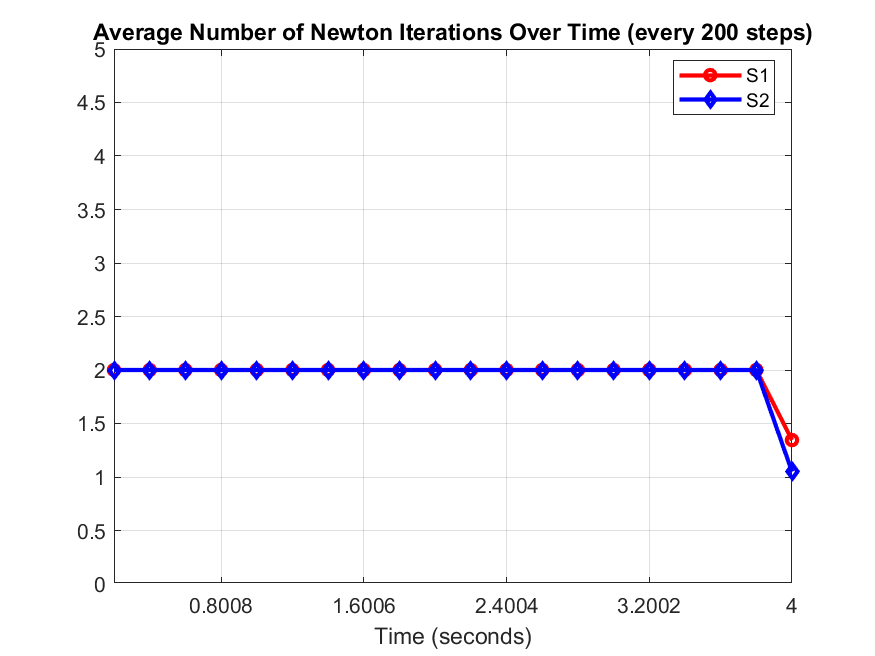}
         \caption{The average number of Newton iterations for schemes S1 and S2, measured every 200 computational steps from $t_0=0$ to $T=4$.}
         \label{fig:Newton_Average_Every_200_Steps_plot}
\end{figure}

For our numerical experiments, we adopt an initial state \(\rho(x,y,1)\) to approximate the above-mentioned steady state, $\rho(x,y,\infty)$. Both schemes, S1 and S2, are employed for this purpose, and the resultant solutions at the final time \(T=4\) are presented in Figure \ref{pmevrand2dinfty2}. For the scheme S2, the decomposed energy functional $E_{2,S2}$ is specifically taken as \(E_{2,S2}(\rho) = \int_{\Omega}\rho(\log \rho - 1) \mathrm{~d}\bx\). We also set $C$ as $10$ to keep $E_{1,S2} + C $ positive, where $E_{1,S2} = E - E_{2,S2}$. As illustrated in Figure \ref{fig:Newton_Average_Every_200_Steps_plot}, S1 and S2 exhibit comparable averages in the number of Newton iterations throughout most of the simulation.

\subsection{PME with large $m$ and a slow drift}
\begin{figure}[tbh]
     \centering
     \begin{subfigure}{0.3\textwidth}
         \centering
         \includegraphics[width=\textwidth]{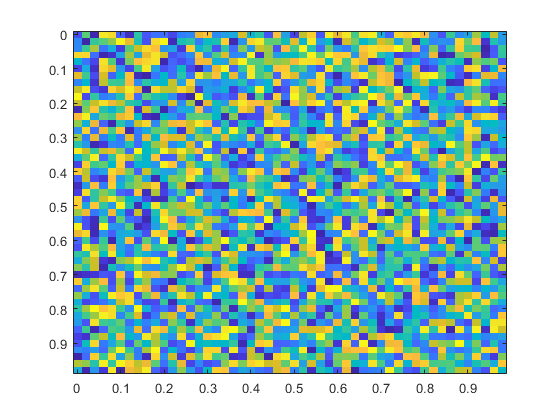}
     \end{subfigure}
     \begin{subfigure}{0.3\textwidth}
         \centering
         \includegraphics[width=\textwidth]{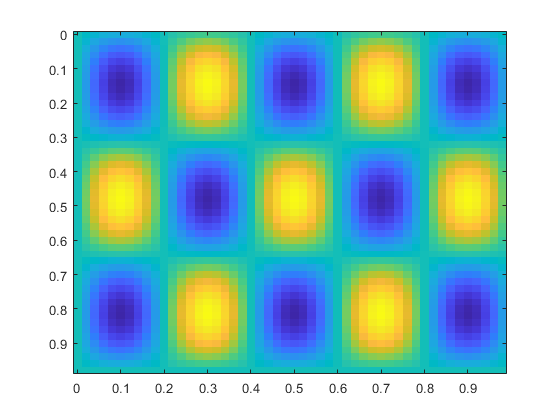}
     \end{subfigure}
    \caption{On the left, the initial state is set to a random configuration using the seed(1). On the right, we depict the potential function given by \(V(x,y) = 1 - \sin(5\pi x)\sin(3\pi y)\).}
     \label{pmevrand2dinit}
\end{figure}

The porous medium equation with a slow drift has earned notable attention in optimal transfers, particularly when accompanied by a high value of \( m \) \cite{jacobs2021back}. A general energy functional for such scenarios can be expressed as:
\begin{equation}
    E(\rho) =  \int_{\Omega} \frac{1}{m-1} \rho^m(\bx) + \rho(\bx) V(\bx) \mathrm{~d}\bx,
\end{equation}
where \( V(\bx) \) is a given potential function. For the scope of this study, we take \( V(x,y) = 1 - \sin(5\pi x)\sin(3\pi y) \). Our simulations are conducted on a domain of \( [-1,1]^2 \), using a grid of \( 50 \times 50 \) spatial points. The initial conditions are randomly set within the \( [-1,1]^2 \) range. The graphical representations of the initial density and the potential function are available in Figure \ref{pmevrand2dinit}.

\begin{figure}[tbh]
     \centering
     \begin{subfigure}{0.32\textwidth}
         \centering
         \includegraphics[width=\textwidth]{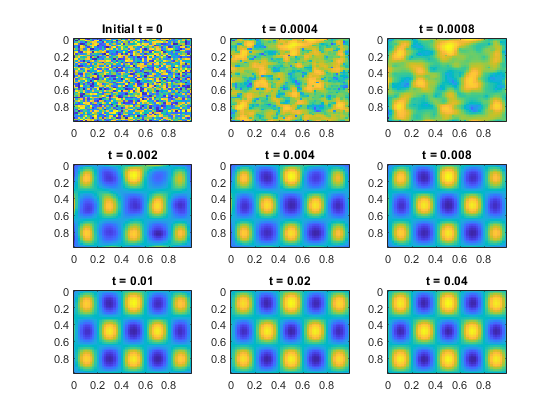}
         \caption{$m=2, T=0.04$}
     \end{subfigure}
     \begin{subfigure}{0.32\textwidth}
         \centering
         \includegraphics[width=\textwidth]{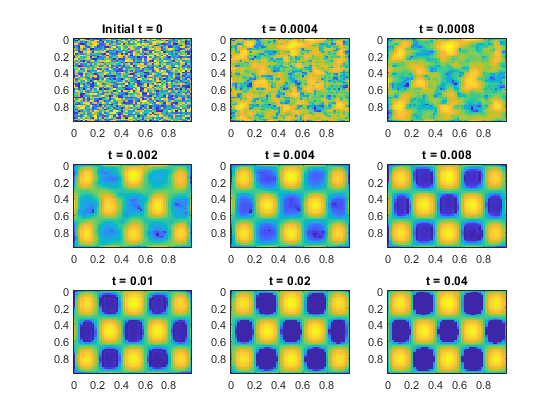}
          \caption{$m=4, T=0.04$}
     \end{subfigure}
     \begin{subfigure}{0.32\textwidth}
         \centering
         \includegraphics[width=\textwidth]{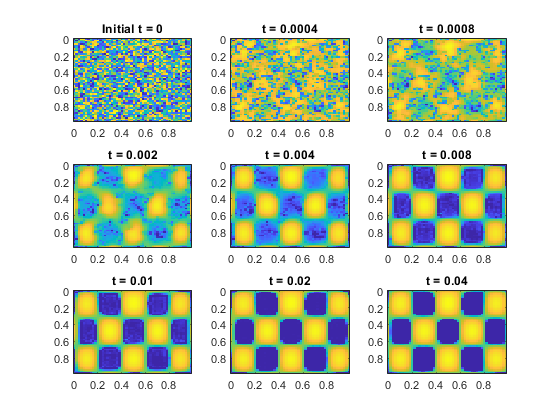}
         \caption{$m=6, T=0.04$}
     \end{subfigure}
     \begin{subfigure}{0.32\textwidth}
         \centering
         \includegraphics[width=\textwidth]{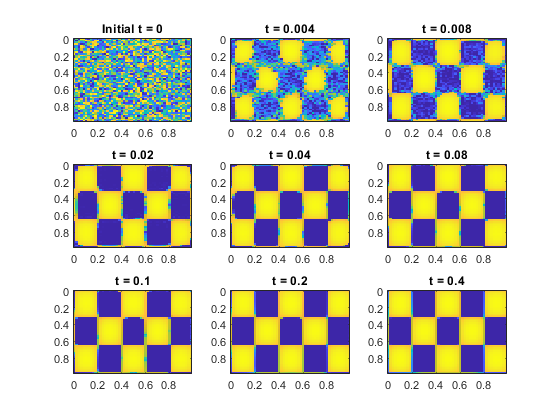}
         \caption{$m=20, T=0.4$}
     \end{subfigure}
     \begin{subfigure}{0.32\textwidth}
         \centering
         \includegraphics[width=\textwidth]{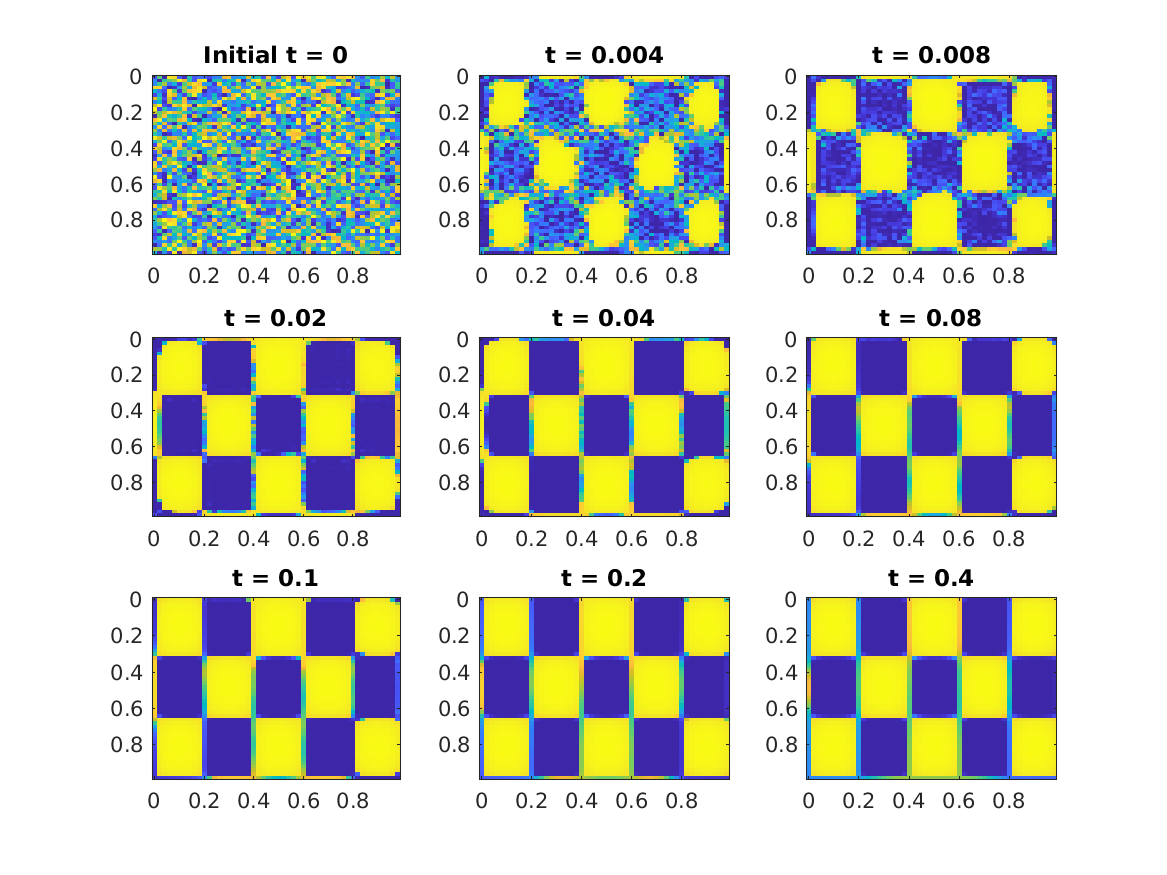}
         \caption{$m=50, T=0.4$}
     \end{subfigure}
     \begin{subfigure}{0.32\textwidth}
         \centering
         \includegraphics[width=\textwidth]{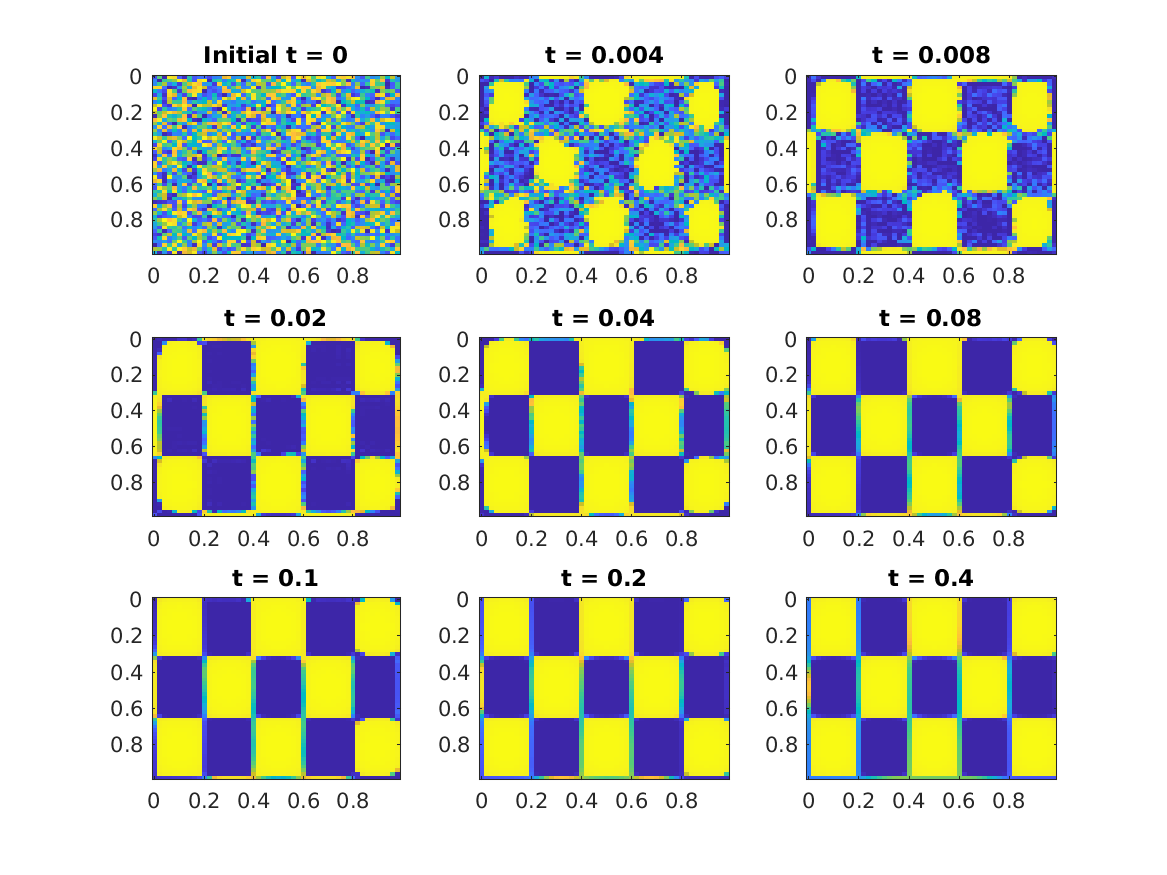}
         \caption{$m=100, T=0.4$}
     \end{subfigure}
     \caption{Evolution of the PME with a slow drift for various values of \(m\). For smaller values of \(m\) (i.e., \(m=2,4,6\)), the time evolution is observed until \(T=0.04\). For larger values (i.e., \(m=20,50,100\)), the evolution is extended up to \(T=0.4\). The results indicate that regions with low potential attract higher density, and the steady state for larger \(m\) values appears more dispersed than with smaller \(m\).}
     \label{pmevrand2dinfty}
\end{figure}

Our scheme S1 does not enjoy the energy dissipation law when $m\geq 2$, whereas the scheme S2 is energy dissipative with respect to a modified energy. Hence, we test our scheme S2 across different values of \( m \): \( m=2, 4, 6, 20, 50, 100\)  to understand how the results change as \( m \) increases. For smaller values of \( m \) (i.e., \( m=2,4,6 \)), we run simulations from \( T=0 \) to \( T=0.04 \) to observe the evolution over time. For larger values (i.e., \( m=20,50,100 \)), we run them from \( T=0 \) to \( T=0.4 \) to see the steady state. The results in Figure \ref{pmevrand2dinfty} indicate that regions with low potential attract higher density and that the steady state with larger $m$ appears to be significantly more dispersed than the steady state with smaller $m$. 

\subsection{Fisher–KPP equation}
As a special case of Onsager gradient flows \eqref{Onsage gf}, we consider the Fisher-KPP equation \cite{bonizzoni2020structure,li2022computational} characterized by potentials \(V_1(\rho)= \alpha \rho\) and \(V_2(\rho)= \frac{\rho(\rho-1)}{2 \log (\rho)}\). The energy associated with this system is given by 
\[
{E}(\rho):=\int_{\Omega}  2 \rho(\log (\rho)-1) \mathrm{~d} \bx + C.
\]
Then the equation can be presented as follows:
\[
\frac{\partial \rho}{\partial t}=\nabla \cdot(2\alpha \rho\nabla \log\rho)+\rho(1-\rho).
\]

We fix the domain \(\Omega = [0,1]\), and set the parameters \(\alpha  = 10^{-4}\), \(C=5\). The initial condition for the density \(\rho\) is defined as \(\rho(x,0) = 0.4\) for \(0 \leq x < 1/2\) and \(\rho(x,0) = 0\) otherwise. For the numerical simulation of the system, we employ the S2 scheme with $E_1 = \int_{\Omega}  \rho(\log (\rho)-1) \mathrm{~d} \bx + C$ and $E_2 = \int_{\Omega}  \rho(\log (\rho)-1) \mathrm{~d} \bx$ combined with a finite difference discretization. In this example, we use \(N = 100\) spatial grid points and a time step size of \(\delta t = 10^{-4}\).

\begin{figure}[tbh]
     \centering
     \begin{subfigure}{0.44\textwidth}
         \centering
         \includegraphics[width=\textwidth]{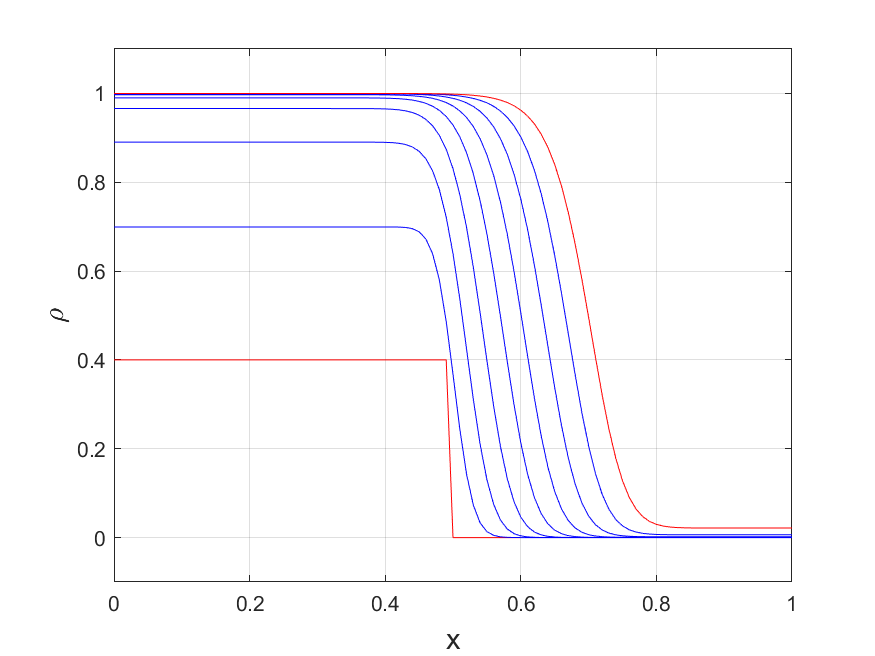}
     \end{subfigure}
     \begin{subfigure}{0.44\textwidth}
         \centering
         \includegraphics[width=\textwidth]{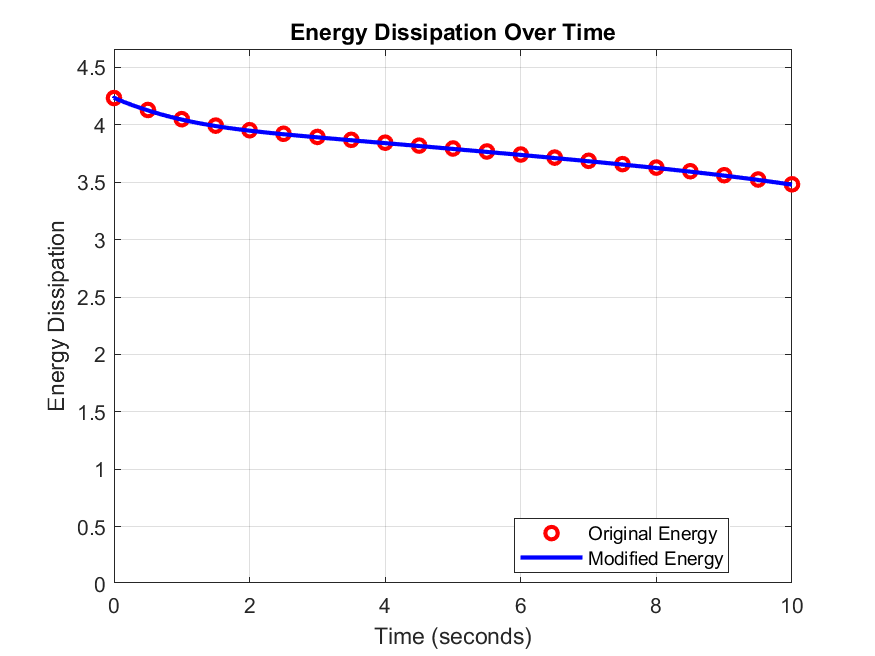}
     \end{subfigure}
    \caption{Evolutionary behavior of the Fisher-KPP equation. The left figure delineates the density evolution from \(t_0=0\) to \(T=10\), emphasizing the initial and terminal states via the red line. The right figure contrasts the original energy and the modified energy dissipation.}
     \label{fkpps3}
\end{figure}

Figure \ref{fkpps3} provides a visual representation of the system's behavior. The left figure illustrates the evolution of the system for time intervals ranging from \(t_0=0\) to \(T=10\), with the initial and final states being demarcated by the red line. The right figure contrasts the modified energy with the original energy, showcasing the energy dissipation as time progresses.

\section{Concluding Remarks}
In this paper, we introduced two novel  numerical schemes specifically tailored for Wasserstein gradient flows. The first scheme is a generalization of the schemes proposed in \cite{ShenXu2021} and \cite{gu2020bound}, while  the second scheme is based on energy splitting along with  a scalar auxiliary variable to ensure energy dissipation.
We demonstrated that both schemes ensure mass conservation, positivity preserving, unique solvability; and the first scheme is  energy dissipative in some special cases while the second scheme is energy dissipative with a modified energy. 

These schemes were designed to address the challenges associated with Wasserstein gradient flows, particularly in preserving positivity and energy dissipation. Each scheme was rigorously tested through a series of numerical experiments, affirming their theoretical precision and computational efficiency. The results confirmed that our schemes not only align with theoretical predictions but also demonstrate significant computational improvements.

In summary, the schemes  proposed in this work are both  robust and practically efficient for solving a class of Wasserstein gradient flows, paving the way for further exploration  in diverse scientific and engineering fields.
\section*{Acknowledgments}
This work is supported in part by NSFC 12371409. This is no conflict of interests.
\bibliographystyle{siamplain}
\bibliography{references}
\end{document}